\documentclass{article}
\usepackage{amsmath,amssymb,graphicx} %load extra symbols and environments
\usepackage{amssymb}
\usepackage{amsthm}
\usepackage[margin=0.9in]{geometry} %set margins

\usepackage{mathtools, nccmath}
\DeclarePairedDelimiter{\nint}\lfloor\rceil
\DeclarePairedDelimiter\floor{\lfloor}{\rfloor}

 \usepackage{diagbox}
\usepackage{url}
\usepackage{algorithm}
\usepackage{algpseudocode}
 \usepackage{tikz}
 \usetikzlibrary{decorations.markings}
 \allowdisplaybreaks
 \usepackage{comment} 
 \usepackage{enumerate}
 \usepackage{appendix}
 \usepackage{soul} %for crossing
\usepackage{color}
\usepackage{lipsum}
\usepackage{lineno}
\begin{document}

%\nocite{*} % this command forces all references in template.bib to be printed in the bibliography

\title{Random walks and moving boundaries: Estimating the penetration of diffusants into dense rubbers}

 \author{Surendra Nepal$^{a}\footnote{corresponding author, email: surendra.nepal@kau.se}$\;, Magnus Ögren$^{b,c}$,  Yosief Wondmagegne$^{a}$, Adrian Muntean$^{a}$\\
$^{a}$ Department of Mathematics and Computer Science, Karlstad University,\\
Universitetsgatan 2, Karlstad, 65188, Sweden\\
$^{b}$ School of Science and Technology, Örebro University\\
SE-701 82, Örebro, Sweden \\
$^{c}$ HMU Research Center, Institute of Emerging Technologies,\\
GR-71004, Heraklion, Greece \\
\date{\today}}
\maketitle
\noindent

\begin{abstract}
For certain  materials science scenarios arising in rubber technology, one-dimensional moving boundary problems (MBPs) with kinetic boundary conditions are capable of unveiling the large-time behavior of the diffusants penetration front, giving a direct estimate on the service life of the material. 
In this paper, we propose  a random walk algorithm able to lead to good numerical approximations of both the concentration profile and the location of the sharp front. 
Essentially, the proposed %discretization 
scheme decouples the target evolution system in two steps: 
%(i) by firstly solving via an explicit Euler method the ordinary differential equation corresponding to the evaluation of the speed of the moving boundary, 
(i) the ordinary differential equation corresponding to the evaluation of the speed of the moving boundary is solved via an explicit Euler method, and 
(ii) the associated diffusion problem is solved by a random walk method.  
To verify the correctness of our random walk algorithm we compare the resulting %obtained random walk 
approximations to results based on a finite element approach with a controlled convergence rate. 
Our numerical experiments recover well penetration depth measurements of an experimental setup targeting dense rubbers. 

	\vskip0.5cm
	\noindent \textit{Keywords:}  Moving boundary problem with a kinetic condition, explicit Euler method,  random walk approximation, finite element approximation
		\vskip0.5cm
	\noindent \textit{MSC 2020 Classification:} 65M75, 65M60, 35R37 \vskip0.2cm

%\noindent \textit {PACS:} 81.05.Lg (Rubber),  74.25.Ha (Penetration depth),  05.40.Fb (Random walks), 02.70.Dh (Finite element analysis)
\noindent \textit {PACS:} 81.05.Lg,  74.25.Ha,  05.40.Fb, 02.70.Dh
\end{abstract}
\section{Introduction}
Rubber is an intensively used material in nowadays technological solutions, e.g., in the context of energy harvesting.
Think for instance of offshore windmills  where rubber joints and ligaments are exposed to diffusants (e.g. chlorides) able to penetrate the porous structure of the material. Even in a dense form, rubbers are to some extent permeable and  many factors affect their permeability. For rubber, permeation is the rate at which small molecules of a gas or liquid transfer through a rubber compound. Although these rates are typically very low, they are important when designing a sealing process. As a general rule, the higher the permeation, the more affected is the service life of the material; see e.g. \cite{borges2021effect, rostami2021chemistry, yasser2023experimental}.

In this framework, we study a random walk method (RWM) to the solution of a one-dimensional moving-boundary problem with a kinetic boundary condition that describes the penetration of a population of diffusants into a dense rubber. The target problem is well-posed (cf. \cite{NHM,kumazaki2020global}), with the solution following foreseeable paths for large times (cf. \cite{ZAMM}).  Furthermore, the approximation by the  finite element method (FEM) of its unique weak solution is well understood  (cf. \cite{nepal2021error, nepal2023analysis}). We did choose to explore the random walk approximation route as an alternative to FEM driven by modeling reasons: It could well be that in many practical situations even if the rubber sample is exposed to an infinite  reservoir of diffusants, these diffusants might find it difficult to penetrate the material and hence only a tiny fraction of them succeed simply to ingress rubber. Having in mind such situations,  one may inquire why a continuum model is needed to describe such ingress. Should such situations appear, then  the random walk approximation becomes the actual model, otherwise approaches like \cite{chester2015finite,wilmers2015continuum} are viable continuum-level alternative descriptions.

The RWM proposed in this framework is rather elementary, its non-standard part is linked to the way we treat the kinetic boundary condition.  We refer the reader to \cite{salsa2016partial} for the basic idea how the method works as well as, for instance, the works \cite{haji1967solution}, \cite{talebi2017study}, \cite{Hogren2014local}, \cite{ogren2020stochastic}, \cite{ogren2019numerical},  \cite{suciu2021global}, \cite{CASABAN2022}, where concrete problems have been successfully solved by RWM and related approaches.  On the computational side of things, it is worth noting that 
for large and complicated physical domains in higher dimensions, the large number of mesh points in a FEM formulation may result in equation systems that are impractical for today’s computers. Instead, random walk methods need only a small RAM memory and can deliver local approximation of solutions at points of interest in the domain.
For domains in practical applications connected to processing digital images, such as photos or CT scans, filtering of the data may result in topological errors that hinder the meshing procedure for FEM. On the contrary, using a random walk method on a large domain with just a few topological errors, will generally result in small statistical errors on the calculated observables.

The paper has the following structure: In Section~\ref{modelequation}, we  describe the governing equations of our model. 
We  present the dimensionless form of the proposed model equations. 
The discussion of the setting of the model equations and the identification of the characteristic scales of dimensionless variables is based on \cite{nepal2021moving}.
In Section~\ref{rwm}, we construct a numerical method that combines the random walk method with the explicit Euler method to solve the dimensionless model equations.   Section~\ref{simulation} contains a couple of numerical experiments and comparisons of random walk solution to the finite element approximation and experimental lab data.   
Our conclusions on the constructed numerical scheme  are listed in Section~\ref{conclusion} together with a couple of suggestions for further investigations.

\section{Setting of the model equations} \label{modelequation}

 For a one-phase Stefan problem with a kinetic condition, we consider the following problem setting. For a fixed given observation time $T_f\in (0, \infty)$, let the interval $[0, T_f]$ be the time span of the involved physical processes.  Let $t \in [0, T_f]$ denote the time variable, $s(t)$ the position of the moving boundary at time $t$ and  $x\in [0, s(t)]$ the space variable, while  $m(t, x)$  is the concentration of diffusant placed in position $x$ at time  $t$.
The diffusants concentration  $m(t, x)$  acts in the non-cylindrical parabolic domain $Q_s(T_{f})$ defined by $$ Q_s(T_{f}):= \{ (t, x) | t \in (0, T_{f}) \; \text{and}\; x \in (0, s(t))\}.$$ 
The problem is: Find 
the profile of 
$m(t, x)$ simultaneously with the position of the moving boundary $x = s(t)$, where the couple $(m(t, x), s(t))$ satisfies the following system of equations:
\begin{align}
\label{4a11}&\displaystyle \frac{\partial m}{\partial t} -D \frac{\partial^2 m}{\partial x^2} = 0\;\;\; \ \text{in}\;\;\; Q_s(T_{f}),\\
\label{4a12}
&-D \frac{\partial m}{\partial x}(t, 0) = \beta(b(t) -\text{H}m(t, 0))  \;\;\; \text{for}\;\;  t\in(0, T_{f}),\\
\label{4a13} &-D \frac{\partial m}{\partial x}(t, s(t))  =s^{\prime}(t)m(t, s(t))  \;\;\; \text{for}\;\; t\in(0, T_{f}),\\
\label{4a14}&s^{\prime}(t) = a_0 (m(t, s(t)) - \sigma(s(t)) )\;\;\;\;\text{for} \;\;\; t \in (0,T_{f}),\\
\label{4a15}&m(0, x) = m_0(x) \;\;\;\text{for}\;\;\; x \in [0, s(0)],\\
\label{4a16} & s(0) = s_0>0\; \text{with}\;\; 0<  s_0< s(t) < \ell.
\end{align}
Here  $a_0>0$ is a  kinetic coefficient, $\beta$ is a positive constant describing the resistance of the $x=0$ interface with respect to the ingress of the diffusants,  $D>0$ is an averaged diffusion constant, and $\text{H}>0$ is the Henry constant. Additionally, $\sigma$ is a function on $\mathbb{R}$, $b$ is a given boundary diffusant concentration on $[0, T_f]$,
$s_0$ is the initial position of the moving boundary, while $m_0$ represents the concentration of the diffusant at $t=0$. As an essential part of this type of modeling, the function $\sigma$ incorporates in an {\em ad hoc} way eventual hyper-elastic effects, i.e. without modeling explicitly the balance of linear momentum for the target material.  

By introducing the dimensionless quantities  
 $z = x/x_{ref}, \hspace{0.4mm}  \tau = tD/x_{ref}^2, \hspace{0.4mm} u = m/m_{ref}, \hspace{0.4mm}   h = s/x_{ref}$, we can define the region for the dimensionless diffusant concentration  by
 $$ Q_h(T):= \left\{ (\tau, z)\, |\, \tau \in (0, T) \; \text{and}\; z \in (0, h(\tau)) \;\text{with}\; T = T_fD/x_{ref}^2 \right\},$$
where $x_{ref}$ is a characteristic  length scale,  while $m_{ref}$ is  a reference value for the diffusant concentration.
In dimensionless form, we can write \eqref{4a11}--\eqref{4a16} as follows: 
\begin{align}
 \label{4a17}&\frac{\partial u}{\partial \tau}  - \frac{\partial ^2 u}{\partial z^2} = 0 \;\;\; \ \text{in}\;\;\; Q_h(T),\\
 \label{4a18}& -\frac{\partial u}{\partial z}(\tau, 0) = \text{Bi} \left(\frac{b(\tau)}{m_{ref}} - {\rm H}u(\tau, 0)\right) \;\;\; \text{for}\;\; \tau\in(0, T), \\
\label{4a19}&-\frac{\partial u}{\partial z}\left(\tau, h(\tau)\right) =h^{\prime}(\tau) u\left(\tau, h(\tau)\right)   \;\;\; \text{for}\;\; \tau\in(0, T),\\
\label{4a20}&h^{\prime}(\tau) = A_0 \left( u \left(\tau, h(\tau)\right) - \frac{\sigma(h(\tau))}{m_{ref}}\right)   \;\;\; \text{for}\;\; \tau\in(0, T),\\
\label{4a21}&u(0, z) = u_0(z) \;\;\; \text{for}\;\;\; z\in[0,  h(0)],\\
\label{4a22}& h(0) = h_0>0\;\; \text{with} \;\;0<h_0<h(\tau)< L ,
\end{align}
where  $u_0(z):=m_0/m_{ref}$, $h_0 := s_0/x_{ref}$ and $L:=\ell/x_{ref}$.  It is worth mentioning that 
 $\text{Bi} := \beta x_{ref}/D$ is the standard  mass transfer Biot number while 
$ A_0 :=  x_{ref} m_{ref} a_0/D$ is the Thiele modulus (or the 2nd Damk\"ohler number). We refer the reader to \cite{nepal2021moving} for more details on the derivation of dimensionless model equations and finite element simulation results.  The convergence analyses of the semi-discrete Galerkin scheme in space and of the fully discrete Galerkin-Euler scheme for the model equations are reported in \cite{nepal2021error, nepal2023analysis}.  

\section{Random walk method}\label{rwm}
In this section, we study a symmetric random walk method (RWM) to solve \eqref{4a17}--\eqref{4a22}, with a suitable bias when treating \eqref{4a20}--\eqref{4a22}. In this method, a finite number of mass-carrying walkers are utilized  to replace the continuous diffusant concentration field.  
The random movement of the walkers ensures that the overall behavior resembles that of macroscopic diffusion. By keeping track of the detailed movements of the random walkers, where they start from, and where they end up, one can understand the dynamics of the problem. We refer the reader to \cite{salsa2016partial} for more details on the theory of  symmetric random walks. Note though that our situation is somewhat special as the walkers positioned at the location of the moving boundary are exposed to a biased random walk; we will explain this particular feature later on in Section \ref{ode_discrete} and Section \ref{interface_discrete}. 

We discretize the space and time domains in the following way:
Let $N, M\in \mathbb{N}$ be 
%a given integer. 
given. We divide the interval $[0, L]$  into $N$  subintervals.
We set $0 =:z_0 < z_1 < \cdots < z_{N} := L$ as discretization points. We define $\Delta z := z_{i+1} - z_{i}$  for $i \in \{0,1,\cdots, N-1\}$ to be a uniform space mesh size. 
We consider a finite sequence of time nodes $0=:\tau_1 < \tau_2 < \cdots < \tau_M := T$ with the uniform time step $\Delta \tau := \tau_{j+1} - \tau_{j}$ for $j \in \{0,1,\cdots, M-1\}.$ 
 We consider walkers walking randomly  along the $z$ axis starting at $z = z_0$. The diffusant concentration $u$  is represented by a discrete number of walkers in the domain. Each walker has a specified mass. For the sake of simplicity, and without loss of generality, we  assume here that all walkers have exactly the same mass of the diffusants concentration of $1$ unit in the box of length $\Delta z$.
The rule is that  each walker draws a uniformly distributed random number and according to the drawn random number decides to go to the new position on the lattice of discrete space nodes. When walkers reach a new position, the new position  gets an increase of $1$ unit of diffusants concentration. More clearly, at each time
step  $\Delta \tau$,  we consider a  unit mass of  diffusants concentration that moves one step of $\Delta z$ unit length randomly to the left or to the right equally likely, that is both with a probability $P = 1/2$. 

Given the initial distribution of the walkers at $\tau = 0$,
we are interested in finding  
the number of walkers starting at  $z = z_0$ for $\tau>0$ and reaching the point $ z = z_i$ after $M$ moves.
To be more precise, we denote the number of walkers at  $(\tau, z) = (\tau_j, z_i)$ by $N_i^j$. 
At the next time 
 step $\tau = \tau_{j+1}$, we then expect $PN_i^j$ walkers moving to the left $ z = z_{i-1}$, and the same number of walkers moving to the right $ z = z_{i+1}$. 
At the same time, $PN_{i-1}^j$ walkers move from the left to $(\tau_{j+1}, z_i)$ and $PN_{i+1}^j$ walkers move from the right to $(\tau_{j+1}, z_i)$. We can express this scenario in terms of the following discrete balance equation:
\begin{align*}
N_i^{j+1} = N_i^{j} - PN_i^j - PN_i^j + PN_{i+1}^j + PN_{i-1}^j.
\end{align*}
 We rearrange the terms to get 
 \begin{align}
\label{4a24}\frac{N_i^{j+1} - N_i^{j}}{\Delta \tau} = d \frac{( N_{i-1}^j- 2N_i^j + N_{i+1}^j)}{(\Delta z)^2}, 
\end{align}
where $d := (P/ \Delta \tau) (\Delta z)^2$. We can see \eqref{4a24} as a  discrete form of  the diffusion equation for $N(\tau, z)$ with the diffusion coefficient $d$. Comparing  diffusion coefficients in \eqref{4a24} and \eqref{4a17} gives the relation 
\begin{align}
  \label{4a25}  \Delta z := \sqrt{2 \Delta \tau}.
\end{align}
 We have described so far the random walk method for a setting where walkers are free to move on an infinite domain. However, we are in fact interested in confining their motion in between a fixed and a moving boundary as described by our original problem posed at the continuum level.  Consequently, our goal is now to solve for $u$ and $h$ in the moving domain where  the moving boundary is driven by \eqref{4a20}.  To be able to solve \eqref{4a17}--\eqref{4a22},  we first decouple the equation for the diffusant concentration  and for the position of the moving boundary. We  then solve them separately. The first step is to approximate the solution of the ordinary differential equation \eqref{4a20} by using the explicit Euler method. The second step is to solve the diffusion equation by the random walk method. The details are presented in algorithm \ref{alg:cap}
 
 Let $h_j = h(\tau_j)$, i.e.   the random walk approximation of $h$ at $\tau_j$ and $u(\tau_j, z_i)$ be the random walk approximation of $u$ at $\tau=\tau_j$ and $z=z_i$. Before continuing with the description of the numerical algorithm,  it is worth defining an index $k_j$      by 
\begin{align}
\label{index}
k_j := \floor*{\frac{h_{j}}{\Delta z}}, \;\,\, j \in \{0, 1,  \cdots, M-1 \},
\end{align}
where $\floor*{x}$ rounds $x$ down towards the nearest integer.
By introducing \eqref{index}, one can keep track of the position of a moving boundary and identify the neighbour discrete  space node within the domain, which can be advantageous for implementing initial and boundary conditions.
At each time slice $\tau = \tau_j$, our system thus consists of a set of discrete space nodes $z = z_i$ for $i \in \{0, 1, \cdots, k_j\}$ in which walkers are distributed.  
For the walker arrived at $z = z_i$ for $i \in \{1, 2, \cdots, k_j-1 \}$ at $\tau = \tau_j$, a random step, based on the value of $p \in \{-1, 1\}$ drawn by each walker with equal probability, is added to the position of all the walkers. 
We then update the number of walkers at each discrete space node for the next time step  by using the following formula:
\begin{align}
u(\tau_{j+1}, z_{i+p}) = u(\tau_{j+1}, z_{i+p}) + 1.
\end{align} 
Depending on the specified boundary condition, a set of rules for the walkers needs to be defined when they reach the boundaries $z = z_0$ and $z= z_{k_j}$.

 \subsection{Treatment of the initial condition}\label{ic}
For the first time step, the initial distribution of walkers within the boundaries is sampled from the initial condition. Otherwise, it is taken from the distribution available from the previous time step.
What concerns our problem, the initial distribution of the walkers is determined by the initial data $u(\tau_0, z_i) = u_0(z_i)$ for any $i\in \{0, 1,  \dots, k_0\}$. To get improved numerical results, we make use of the following scaling argument: We multiply the number of walkers starting at all discrete points
defined by the initial value with a large integer $n$, and finally, we divide the diffusants concentration profile by $n$. That means the total number of walkers introduced at $\tau = 0$ is given by $n u_0(z_i)$ for   $i \in \{0, 1, \dots, k_0\}$.  
\subsection{Treatment of  Dirichlet and Robin boundary conditions at the fixed boundary} \label{fixedbc}
To approximate numerically  the non-homogeneous Dirichlet boundary condition on the left boundary, we employ the following two-step procedure:
 The Step 1 involves removing or absorbing all walkers that reach the left boundary, which is equivalent to a homogeneous Dirichlet condition. In Step 2, we adjust the number of walkers at the left boundary during each time step to meet the non-homogeneous boundary condition.  The walkers that remain inside the boundary after a time step are allowed to continue their random walk during the next time step, while those who cross the left boundary are eliminated. The described two-steps procedure  is sometimes referred to as  the Söderholm procedure; this is a possible strategy  for implementing  the non-homogeneous Dirichlet boundary condition in the random walk method. We refer the reader to \cite{schwind2003some, schwind2001random} for more details.
 
To deal with a Robin boundary condition at the left boundary (an inflow boundary condition), we introduce walkers starting at the left fixed boundary.   
The number of walkers starting at the left boundary $u(\tau_j, z_0)$  is computed based on \eqref{4a18}. We use  the following forward difference approximation for approximating the space derivative of $u$, 
\begin{align*}
\frac{\partial u}{\partial z}(\tau_j, z_0) \approx \frac{u(\tau_j, z_1)- u(\tau_j, z_0) }{\Delta z}.
\end{align*}
From \eqref{4a18}, we can write 
\begin{align*}
\frac{u(\tau_j, z_0)- u(\tau_j, z_1) }{\Delta z} = \text{Bi} \left( \frac{b(\tau_j)}{m_0} - \text{H}u(\tau_j, z_0) \right).
\end{align*}
We multiply the number of walkers  defined by   the boundary conditions with the number $n$. 
We then update the number of walkers at the left boundary according to the rule:
\begin{align}
\label{4a28}
u(\tau_j, z_0) = \nint*{\frac{n\Delta z \text{Bi}\, b(\tau_j)/ m_0 +  u(\tau_j, z_1)}{ 1+ \Delta z \rm{Bi}\, \text{H}}} \;\;\; \text{for}\; j \in \{1,2, \cdots, M\},
\end{align}
where $\nint{x}$ rounds $x$ to the nearest integer. 
 
\subsection{Treatment of the ordinary differential equation}\label{ode_discrete}
   
 A key issue in solving moving boundary problems by the random walk method is the handling of the moving boundary position and its velocity (or its speed for the 1D case). Here, the movement of the boundary $ h(\tau)$  needs to be traced together with the approximation of the concentration of diffusants. We discretize the ordinary differential equation  defining the speed of the moving boundary by using the explicit Euler method.
 We then obtain the  random walk approximation of $h$ at $\tau_{j+1}$  as 
\begin{align} \label{4a26}
 h_{j + 1} = h_{j}  + \frac{\Delta h_j}{n}, \;\,\, j \in \{0, 1, 2, \cdots, M-1 \},
\end{align}
 where the total increment of the position of the moving boundary  for each time step  is defined by 
\begin{align}
\label{4a29}\frac{\Delta h_j}{n} := N_j\left[\frac{\Delta  \tau A_0}{n} \left( u \left(\tau_j, h_j\right) - \frac{\sigma(h_j)}{m_{ref}}\right) \right]. 
\end{align}
Here $N_j$ is the total number of walkers who arrived at the moving boundary at the time $\tau = \tau_j$ and contribute to the increment of the moving boundary. The term 
\begin{align} 
\frac{\Delta  \tau A_0}{n} \left( u \left(\tau_j, h_j\right) - \frac{\sigma(h_j)}{m_{ref}}\right),
\label{eq:IncrementOfTheBoundary}
\end{align}
 is the increment  of the boundary for a walker at time $\tau = \tau_j$.  In \eqref{4a26},  we  divide by $n$ to adjust for the multiplication with the factor $n$ at the starting points, as discussed in Section \ref{ic} and  Section \ref{fixedbc}. It is worth mentioning that  $N_j$ is unknown at this moment. But, it will be identified while discussing the Robin boundary condition in  Section \ref{interface_discrete}.  Once we know the value for $N_0$,  we can compute $h_1$ as $\Delta h_0$  is computable from given initial conditions. We then solve the diffusion equation in the fixed domain $(0, T) \times (0, h_1)$ by the random walk method. We  compute $\Delta h_1$ and $h_2$ and then solve for $u$. We continue the same process until the final time $T$ is reached or until all the walkers arrived at the boundary $L$.

The increment in the moving boundary $\Delta h_j$ varies  over time and may not be equal to the space mesh size $\Delta z$. Hence,  at any time $\tau_j = j\Delta \tau$  the current position of the moving boundary may not coincide with a space mesh point. 
To get to know the distance between the moving boundary position  and its neighbouring discrete space point inside the moving domain, we defined the index $k_j$ for $ j \in \{0, 1, \cdots, M-1 \}$  in  \eqref{index}. 
The index \eqref{index} defines  the walker's decision to move or not towards the current location of the moving boundary.  If the distance between the position of boundary $h_{j}$ and the position of the walker  is less than $\Delta z$, then the walker does not have a place to move to the right. That is if the inequality
$0\leq h_{j} - z_{k_j} \leq \Delta z$ holds true, then the walker moves to the left, or  the boundary progresses according to \eqref{4a29},  depending on the boundary condition at the moving boundary.  While implementing \eqref{4a29}, we use $u\left(\tau_j, z_{k_j}\right)$ instead of $u \left(\tau_j, h_j\right)$. 

\subsection{Treatment of the moving boundary condition}\label{interface_discrete}

At the right boundary $z= h(\tau)$, we have a Robin-type boundary condition \eqref{4a19}. 
This type of boundary condition is suitable for a partial reflection scenario; see e.g. \cite{Hogren2014local}. In other words, 
 a fraction of the walkers hitting the boundary is adsorbed with some reaction probability $P_b$, while the rest of the  walkers are reflected to the left. 
 
 When the walkers succeed to arrive close to the moving boundary, i.e. at $z = z_{k_j}$, then the following  rule applies:
\begin{enumerate}[Step 1:]%[label=  Step \arabic*:]
    \item The walker present at $ z = z_{k_j}$ moves to the left if  $p = -1$. If $p=1$, then Step 2 and Step 3 below are followed.     
    \item For all walkers arrived at node $z = z_{k_j}$ and $p=1$, we evaluate the reaction probability $P_b(\tau_j)$ at that node and compare it to a uniformly generated random number $r$ between $(0,1)$.  If $r$ is less than the current evaluation of the reaction probability  $P_b(\tau_j)$, then the walkers  can potentially change  the position of the moving boundary. The advancement of the boundary continues unless the moving boundary crosses the right-hand side neighboring the discrete space node $z = z_{k_j+1}$. Once the moving boundary crosses this node, the new discrete node $z = z_{k_j+1}$ inside the domain is updated as a node close to the moving boundary.  Consequently, the remaining walkers at $z = z_{k_j}$ move  either to the left  $z = z_{k_j-1}$ or to the right $z = z_{k_j+1}$ with equal probability. This means that only a fraction of the walkers arrived at $z = z_{k_j}$, denoted by $N_j$ in \eqref{4a29}, contribute to the overall increase or decrease of the boundary and the rest of them contribute to the diffusion process. 
    \item Walkers that satisfy the condition $r\geq P_b(\tau_j)$ are reflected to the left and those that contribute to the  increment of the moving boundary are adsorbed.  They stick to the boundary.  
\end{enumerate}
In the context of our problem, we make the following splitting: firstly, we allow the walkers to diffuse on the left side from the point close to the moving  boundary, and secondly, the walkers are exposed to the interplay between the kinetic condition and the diffusion process.
This means that the walkers are in principle able to "push" the moving boundary unless  breaking/acceleration  effects, due to the term $\sigma(\cdot)$, dominate the diffusion process. When the cumulative effect of their individual contributions leads to an overall increase and the walker crosses the neighboring discrete space node on the right-hand side,  the diffusant diffuses further up to that neighboring discrete space node.  
In the third step, the walkers decide their next movement based on the values of $r$ and $P_b(\tau_j)$.

The probability that the walkers stay at the same place 
rather than being reflected depends on the reactivity of the boundary $h^{\prime}(\tau_j)$. 
 Relying on the definition of reaction probability defined in \cite{Hogren2014local}, we define the reaction probability $P_b$ for our problem by 
 \begin{align}
 P_b(\tau_j) = \sqrt{2 \Delta \tau} h^\prime(\tau_j).
 \end{align}
 Other choices for the definition of the reaction probability are possible,  we refer the reader to \cite{Erban_2007, boccardo2018improved} for more details on this matter. For each walker  who arrived at the boundary and obtains $p=1$, equation \eqref{4a29} leads to the following formula for computing the probability 
 
\begin{align} \label{4a27}
P_b(\tau_j) = \frac{ \sqrt{2 \Delta \tau} A_0}{n} \left( u \left(\tau_j, z_{k_j}\right) - \frac{\sigma(h_j)}{m_{ref}}\right).
\end{align}
 
We remind the reader that $h_j = h(\tau_j)$. If  $r<P_b(\tau_j)$, then we update  the position of the moving boundary by adding the increment of the moving boundary for a walker~\eqref{eq:IncrementOfTheBoundary}, defined in \eqref{4a29}. 

The number of walkers for the next time step is updated by
\begin{align}
u(\tau_{j +1}, z_{k_j}) = u(\tau_{j +1}, z_{k_j}) +1.
\end{align}
If $r\geq P_b(\tau_j)$, then the walker is reflected and  moves to the left. %For this, we update the 
For a small diffusion coefficient $D$, $P_b(\tau_j)$ can have a large value, see the definition of $A_0$ under \eqref{4a22}.
 Hence, a suitable time step size $\Delta \tau$ must be chosen in order to obtain the value of $P_b(\tau_j)$ between zero and one. \\
Claim: assume  the time step size $ \Delta \tau \in(0,1)$ satisfies \eqref{4a25} and the following condition 
\begin{align} 
\label{4a30}\sqrt{\Delta \tau} \leq  \frac{n}{\sqrt{2} A_0 \max_j \left( u(\tau_j, z_{k_j}) - \frac{\sigma(h_j)}{m_{ref}}\right) }.\end{align}
Then it holds that  $P_b(\tau_j)$ in \eqref{4a27} is positive and bounded above by $1$.
Moreover, the increment of the moving boundary for each walker defined in \eqref{4a29} satisfies the following inequality 
\begin{align}
\label{4a31} 
  \frac{\Delta  \tau A_0}{n} \max_j \left( u \left(\tau_j, z_{k_j}\right) - \frac{\sigma(h_j)}{m_{ref}}\right) <  \Delta z.
\end{align}
 The positivity of $P_b(\cdot)$ depends on the choice of the function $\sigma(\cdot)$. To preserve the positivity, we have to choose $\sigma(\cdot)$ in such a way that the following inequality holds:
 \begin{align*}
 \min_{j}\left( u(\tau_j, z_{k_j}) - \frac{\sigma(h_j)}{m_{ref}}\right)> 0.
 \end{align*}
 In this case, all walkers who satisfy $r<P_b(\tau_j)$ push the boundary. 
In addition, it also holds that $h_{j+1}>h_j$ for all $j \in \{0, 1, \cdots, M-1\}$.  
The upper bound for  $P_b(\cdot)$ (i.e. by $1$)  follows from  \eqref{4a27}. \\ 
Using \eqref{4a30}, we estimate the maximum increment of the moving boundary for each walker by
\begin{align*}
 \frac{\Delta  \tau A_0}{n} \max_j \left( u \left(\tau_j, z_{k_j}\right) - \frac{\sigma(h_j)}{m_{ref}}\right)  \leq \Delta \tau \frac{1}{ \sqrt {2 \Delta \tau}} = \frac{\Delta z}{2} < \Delta z.
\end{align*}
This completes the proof of the claim. 
It is worth mentioning that if the following holds $$\frac{\Delta  \tau A_0}{n}  \left( u \left(\tau_j, z_{k_j}\right) - \frac{\sigma(h_j)}{m_{ref}}\right) \gg  \Delta z,
$$ for some $j$, then the boundary will move several $\Delta z$ steps from the contribution of a walker. This leads to a poor approximation of the position of the  moving boundary. To solve this issue, we have to choose the number $n$ in such a way  that \eqref{4a31} holds true. 
\begin{figure}[ht]
	\begin{center}
		\begin{tikzpicture}[scale=1, every node/.style={scale=1}]
		\draw[step=1cm] (0,0) grid (6,7);
       \draw[very thick,->] (-0.5, 0) -- (6.5,0) node[anchor=north west] {$z$};
       \draw[very thick,->] (0,-0.5) -- (0.0, 7.5) node[anchor=south east] {$\tau$};
       \foreach \x in {1,2,3,4,5,6}
    \draw (\x cm,1pt) -- (\x cm,-1pt) node[anchor=north] {$z_\x$};
  \foreach \y in {1,2,3,4, 5, 6, 7}
    \draw (1pt,\y cm) -- (-1pt,\y cm) node[anchor=east] {$\tau_\y$};
\draw[scale=0.5,domain=2:8.5,smooth,variable=\x,blue,very thick] plot ({\x},{(\x-2)*(\x-2)/3 });
   %\addplot[blue,thick,samples=100] {x^2};
%\addlegendentry{$y=x^2$}
\begin{scope}[very thick,decoration={
    markings,
    mark=at position 0.5 with {\arrow{latex}}}
    ]  
    \foreach \i in {0,...,5, 6}{
    % Draw the grid
    %\draw[gray] (0,\i) -- (7-\i,7);  
    \draw[postaction={decorate},dashed, red] (1,0) -- (2,1);
    \draw[postaction={decorate},dashed, red] (2,2) -- (3,3);
    \draw[postaction={decorate},dashed, red] (3,5) -- (4,6);
    \draw[postaction={decorate},red] (2, 1) -- (1,2);
    \draw[postaction={decorate}, red] (3, 3) -- (2,4);
    \draw[postaction={decorate},red] (0,0) -- (1,1);
    \draw[postaction={decorate},red] (1,1) -- (2,2);
    \draw[postaction={decorate},red] (0,1) -- (1,2);
    \draw[postaction={decorate},red] (1,2) -- (2,3);
    \draw[postaction={decorate},red] (2,3) -- (3,4);
    \draw[postaction={decorate},red] (4, 6) -- (3,7);
     \draw[postaction={decorate},red] (0,2) -- (1,3);
     \draw[postaction={decorate},red] (1,3) -- (2,4);
     \draw[postaction={decorate},red] (2,4) -- (3,5);
     \draw[postaction={decorate},red] (0,3) -- (1,4);
     \draw[postaction={decorate},red] (1,4) -- (2,5);
     \draw[postaction={decorate},red] (2,5) -- (3,6);
     \draw[postaction={decorate},red] (3,6) -- (4,7);
     \draw[postaction={decorate},red] (0,4) -- (1,5);
    \draw[postaction={decorate},red] (1,5) -- (2,6);
    \draw[postaction={decorate},red] (2,6) -- (3,7);
     \draw[postaction={decorate},red] (0,5) -- (1,6);
     \draw[postaction={decorate},red] (1,6) -- (2,7);
     \draw[postaction={decorate},red] (0,6) -- (1,7);
     %%%%%%%%%%%%%%%%%%%%%%%%
     \draw[postaction={decorate},red] (2,2) -- (1,3);
    \draw[postaction={decorate},red] (2,3) -- (1,4);
     \draw[postaction={decorate},red] (2,4) -- (1,5);
     \draw[postaction={decorate},red] (3,4) -- (2,5);
     \draw[postaction={decorate},red] (2,5) -- (1,6);
     %\draw[red] (2,6) -- (1,7);
     \draw[postaction={decorate},red] (3,5) -- (2,6);
     \draw[postaction={decorate},red] (2,6) -- (1,7);
     \draw[postaction={decorate},red] (3, 6) -- (2,7);
    \draw[white] (1, 0) -- (1,1);
      \draw[postaction={decorate},dash dot, red] (1, 0) -- (1,1);
      \draw[white] (2, 1) -- (2,2);
      \draw[postaction={decorate},dash dot, red] (2, 1) -- (2,2);
      \draw[white] (2, 2) -- (2,3);
      \draw[postaction={decorate},dash dot, red] (2, 2) -- (2,3);
      \draw[white] (3, 3) -- (3,6);
      \draw[postaction={decorate},dash dot, red] (3, 3) -- (3,4);
      \draw[postaction={decorate},dash dot, red] (3, 4) -- (3,5);
      \draw[postaction={decorate},dash dot, red] (3, 5) -- (3,6);
      \draw[white] (4, 6) -- (4,7);
      \draw[postaction={decorate},dash dot,red](4, 6) -- (4,7);
    }
    \end{scope}
  \foreach \x in {0}
    \foreach \y in {0,...,7}
    {
    \fill (\x,\y) circle(2pt);
    %\fill (\x+1,\y+2) circle (2pt);
    }
    \foreach \p in {(1, 0), (2,1), (3, 3), (4, 6), (0,0), (0,1),(0,2),(0,3),(0,4), (0,5), (0,6), (0,7), (1, 1), (1,2),(1,3),(1,4), (1,5),(1,6),(1,7),(2,2), (2,3), (2,4), (2, 5), (2, 6), (2,7), (3,4), (3,5), (3,6), (3, 7), (4,7)}
      \fill \p circle(2.5pt);
      \foreach \p in {(1,0), (2,1), (2,2), (3,3),
      (3,4), (3, 5), (4, 6), (4, 7)}
      \node[circle,draw=black, fill=white, inner sep=1pt,minimum size=5pt] at \p {};
		\end{tikzpicture}
		\caption{Sketch of the dynamics of the walkers inside  the moving domain -- a vizualization of  Step~1-- Step~3 in the algorithm of Section~\ref{interface_discrete}.}
		\label{Fig:2}
	\end{center}
\end{figure}
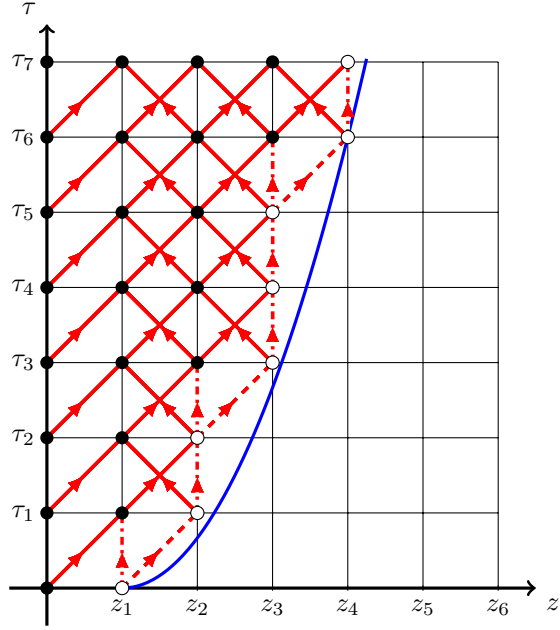

We present  a sketch in Figure \ref{Fig:2} to illustrate the dynamics of the walkers inside the domain. 
The solid line connecting two black-black or black-white nodes represents the movement of the walkers with an equal probability, i.e. $P= 1/2$. 
Walkers who arrived at the white nodes, i.e. $z = z_{k_j}$, follow the moving boundary conditions described in Step 1--Step 3. The solid line connecting two white-black nodes represents the movement of the walkers with a probability greater than or equal to $1/2$, depending on the conditions Step~1 and Step~3.
 The dashed line connecting two white-white nodes represents the movement of the walkers after the moving boundary crosses the neighbourhood node. 
The dashed-dotted line connecting two white-white nodes indicates the movement of the walkers who push the moving boundary.
%\begin{tikzpicture}[scale=3]
%\draw[step=.5cm,gray,very thin] (-2.0,-2.0) grid (2.0,2.0);
%\end{tikzpicture}
%\end{comment}
\begin{algorithm}[ht]
\caption{Procedure to compute the random walk solution to  \eqref{4a17}--\eqref{4a22}.} \label{alg:cap}
\begin{algorithmic}[1]
\State Choose data $a_0, m_0, s_0, D, b, \sigma$.
\State Select a uniform spatial step-size $\Delta z$ generating the nodes $z_i = i \Delta z, 0 \leq i \leq N$, in $[0, L]$ such that $N\Delta z = L$. $N$ is the total number of discrete space nodes in the domain $[0, L]$. 
\State Choose a time step  $\Delta \tau$ satisfying the condition \eqref{4a25}. 
\State Consider a partition of the temporal interval $[0, T]$  generating the nodes $\tau_j = j \Delta \tau, 0 \leq j \leq M$, in $[0, T]$ such that $M\Delta \tau = T$.  $M$ is the total number of discrete space nodes in the domain $[0, T]$. 
\State Initialize a matrix $u$ for diffusants  and a vector $h$ for the moving boundary.
\State Initialize the index for space $i$, time $j$ and  moving boundary $k_j$. 
\While{$ j < M-1$ and $ k_j < N -1$ }
\If{$j == 0$}
    \State Set the given initial value for $u$, i.e., $u[0, 0] = u_0 n$. 
\Else
    \State Compute the number of walker at the left boundary $z= z_0$ by using \eqref{4a28}.  
\EndIf
\ForAll{$i \in \{0, 1, \cdots, N\}$}
\ForAll{$k \in range(u[j, i])$}
\State Generate a number randomly $p = -1$ or $ p = +1$
\If{$i+p > 0$ and $i + p <= k_j$ and $i <= N-1$}
\State Update $u$ i.e., $u[j +1, i+p] = u[j +1, i+p] + 1$  
\ElsIf{$i + p == k_j +1$}
\State Compute $\displaystyle\frac{\Delta h_j}{n}$ using \eqref{4a29}
\State Compute $P_b$ using \eqref{4a27}. 
\State Generate a random number $r$ uniformly distributed in the interval $(0, 1)$.
\If{$r < P_b$}
\State Update $h_{j+1}$ using \eqref{4a26}
\State Walkers stay at the boundary, i.e., $u[j+1, k_j] = u[j+1, k_j] + 1$.
\Else
\State  Walkers move one step to the left, i.e., $u[j+1, k_j-1] = u[j+1, k_j-1] + 1$
\EndIf
\State Update $k_j$, i.e., $k_j = \text{floor}(h_j/\Delta z)$.
\EndIf
\EndFor
\EndFor
\State Update time index, i.e., $j = j + 1$
\EndWhile 
\State Transfer dimensionless form of $u$ and $h$ to dimensional form $m$ and $s$.
\end{algorithmic}
\end{algorithm}

\section{Numerical results} \label{simulation}
 Our goal in this section is to  present the simulation results obtained by our  random walk algorithm.  
 Firstly, we  solve  the  dimensionless model equations  with a non-homogeneous Dirichlet boundary condition at $z=0$ using the random walk method. As a next step,  we  solve the problem presented in \eqref{4a17}--\eqref{4a22}  and exhibit the results of the simulation. Lastly, we solve the dimensional form of the equation described in  \eqref{4a11}--\eqref{4a16} to represent the laboratory experimental data.   An analytical solution to the problem is not available. Thus, we test the performance of the random walk algorithm  in each case by  solving the same problem using the finite element method and comparing the resulting random walk solution with the finite element solution.  To approximate with finite elements,  we first transform the moving boundary problem  to a problem posed in a fixed  domain  using the transformation $y = z/ h(\tau) \in (0, 1)$. We then solve the transformed equations using the finite element methods. For details on the finite element approximation and simulation of \eqref{4a11}--\eqref{4a16}, we refer the reader to our previous work \cite{nepal2021error}.

\subsection{Simulation results for non-homogeneous Dirichlet boundary condition} \label{nonhomdirichlet}
 In this section, we present the random walk solution to \eqref{4a17}--\eqref{4a22} except for the boundary condition at $z=0$. Instead of a Robin-type boundary condition  defined in \eqref{4a18},  we work with  a non-homogeneous Dirichlet boundary condition  defined by 
$$ u(\tau, 0) = u_D  \;\;\; \text{for}\;\; \tau\in(0, T),$$
where $u_D$ is a constant. 
The aim of treating this case first is to make the problem simpler and understand the qualitative behaviour of the random walk solution.  The observation time is set to $T = 0.0001$.   We take the value $10$ mm and $0.5$ gram/mm$^3$ for the characteristic length scale $x_{ref}$ and for reference diffusant concentration $m_{ref}$ respectively.  The function $\sigma(h(\tau))$ is chosen as  $ \frac{h(\tau)}{20}$. This linear  choice for $\sigma(h(\tau))$ is taken from our previous work \cite{nepal2021moving}.   With choosing $h_0 = 0.001$, we compute the index  $k_{0} := \floor*{\frac{h_0}{\Delta z}}$.   Initially, we place $n$ walkers at each discrete node  $z = z_i$ for $ i \in \{0, 1, \cdots, k_0\}$. We choose $u_D = 10$ so that at each time step, we introduce $10n$ walkers at the left boundary $z= 0$. With our choice of parameters, the dimensionless number $A_0$ defined in Section \ref{modelequation} becomes $2500$. 
 
\begin{figure}[ht] 
\centering 
\includegraphics[width=0.45\textwidth]{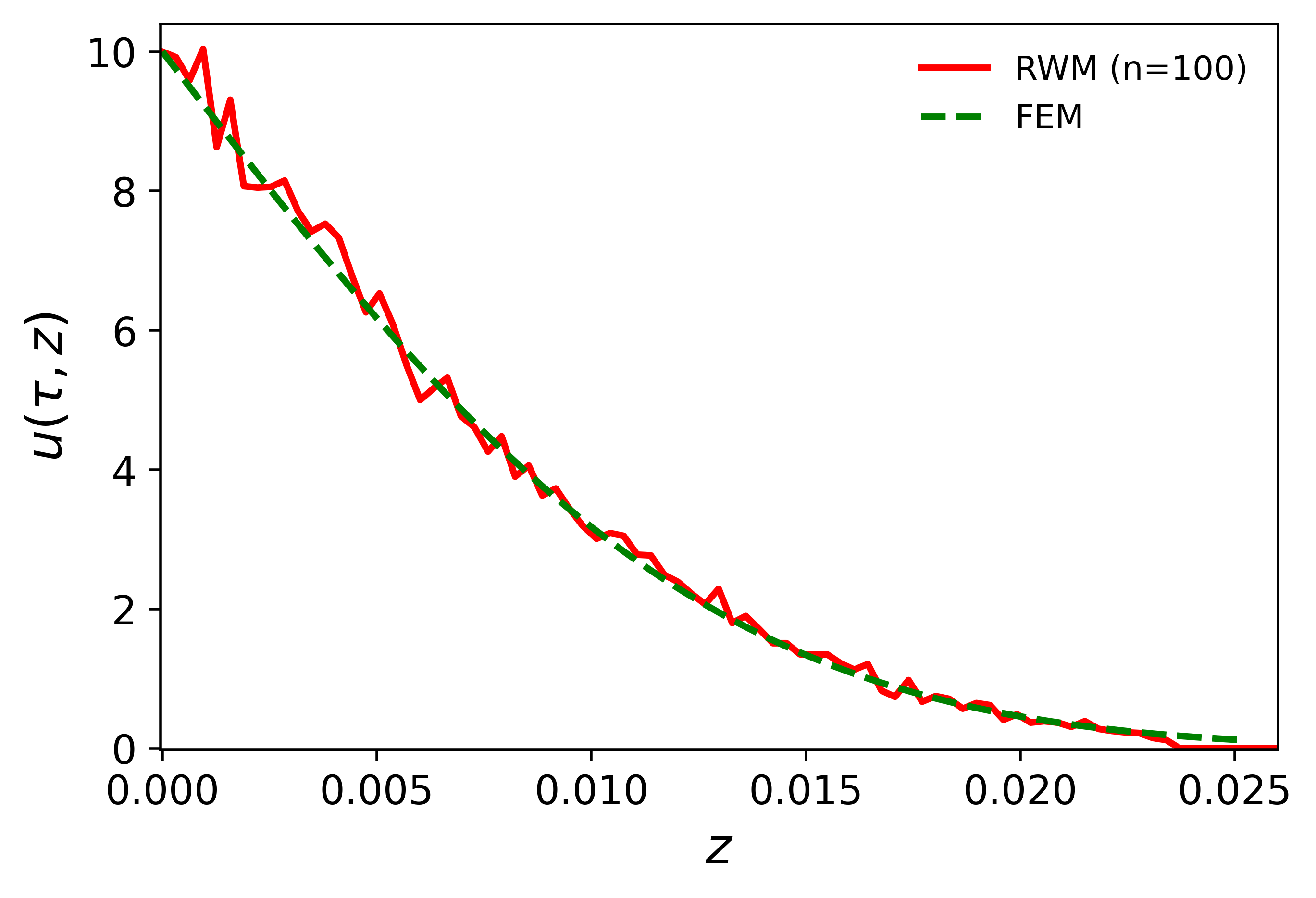}
\hspace{0.01cm}
\includegraphics[width = 0.45\textwidth]{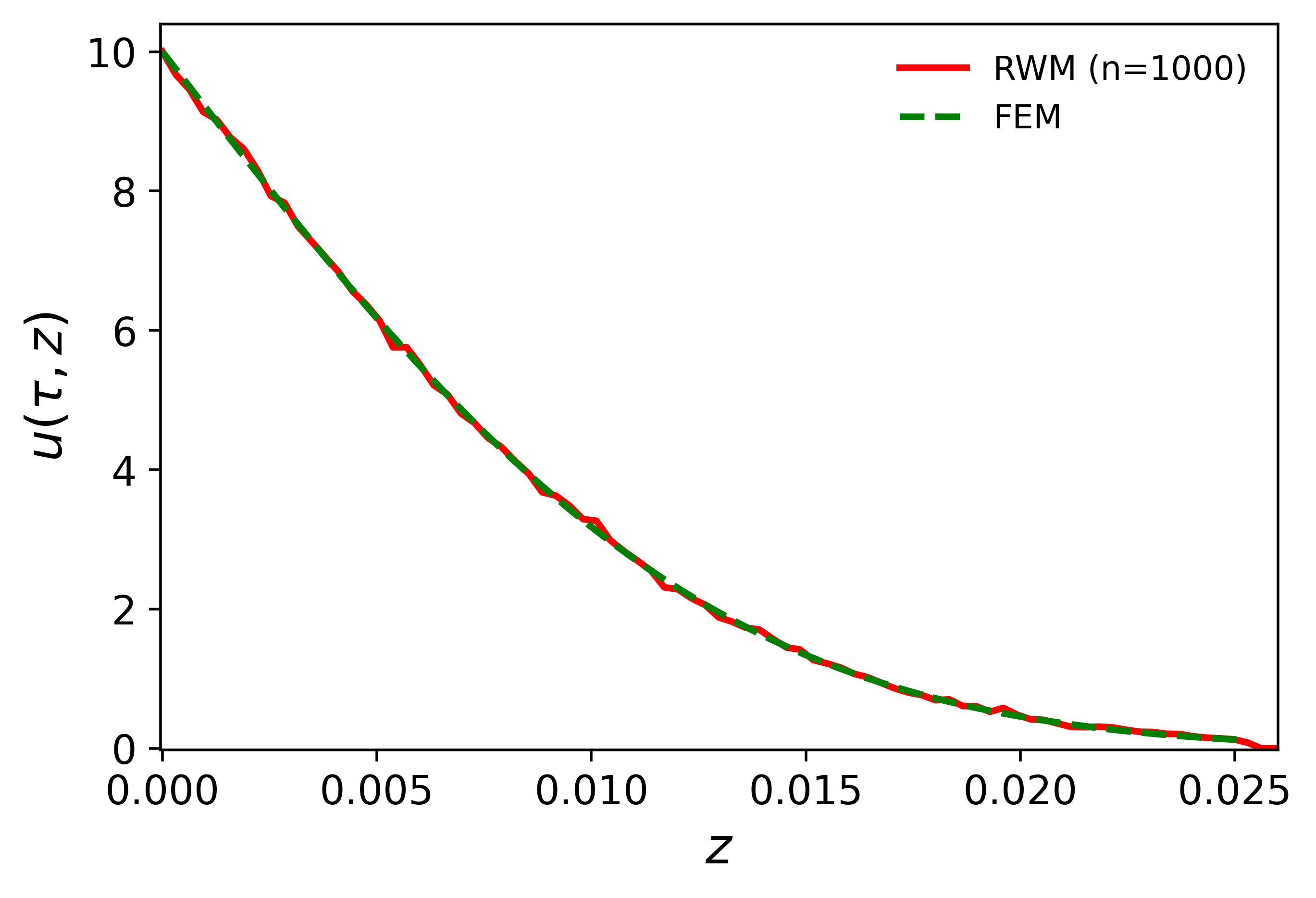}
	\caption{ Concentration profile at $\tau = 0.00005$ by FEM and RWM for different values of $n$, for $n=100$ (left), and $n=1000$ (right),  with  $\Delta \tau = 5 \times 10^{-8}$. }
	\label{compare22}
\end{figure}

\begin{figure}[ht] 
\centering 
\includegraphics[width=0.45\textwidth]{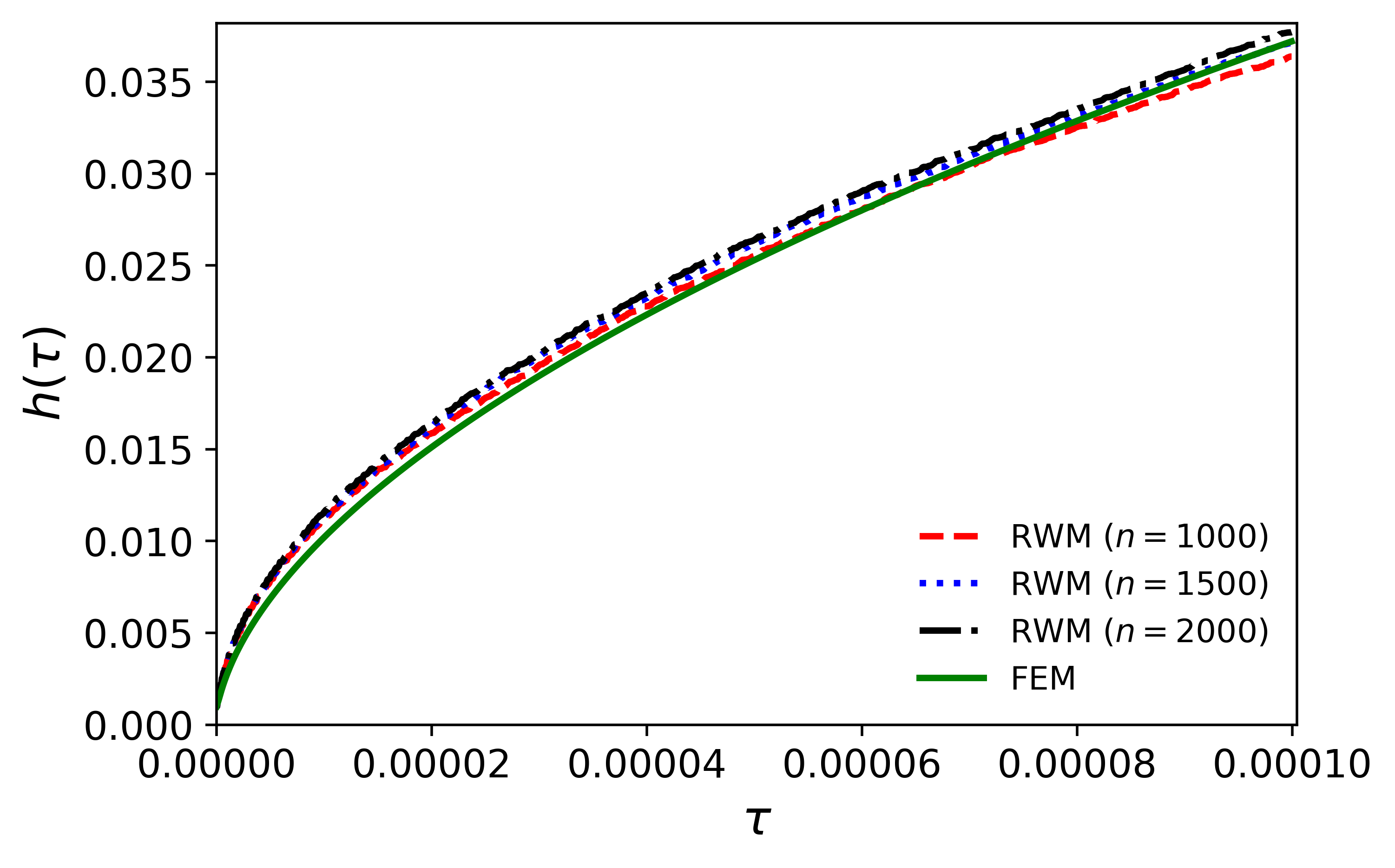}
\hspace{0.001cm}
\includegraphics[width = 0.45\textwidth]{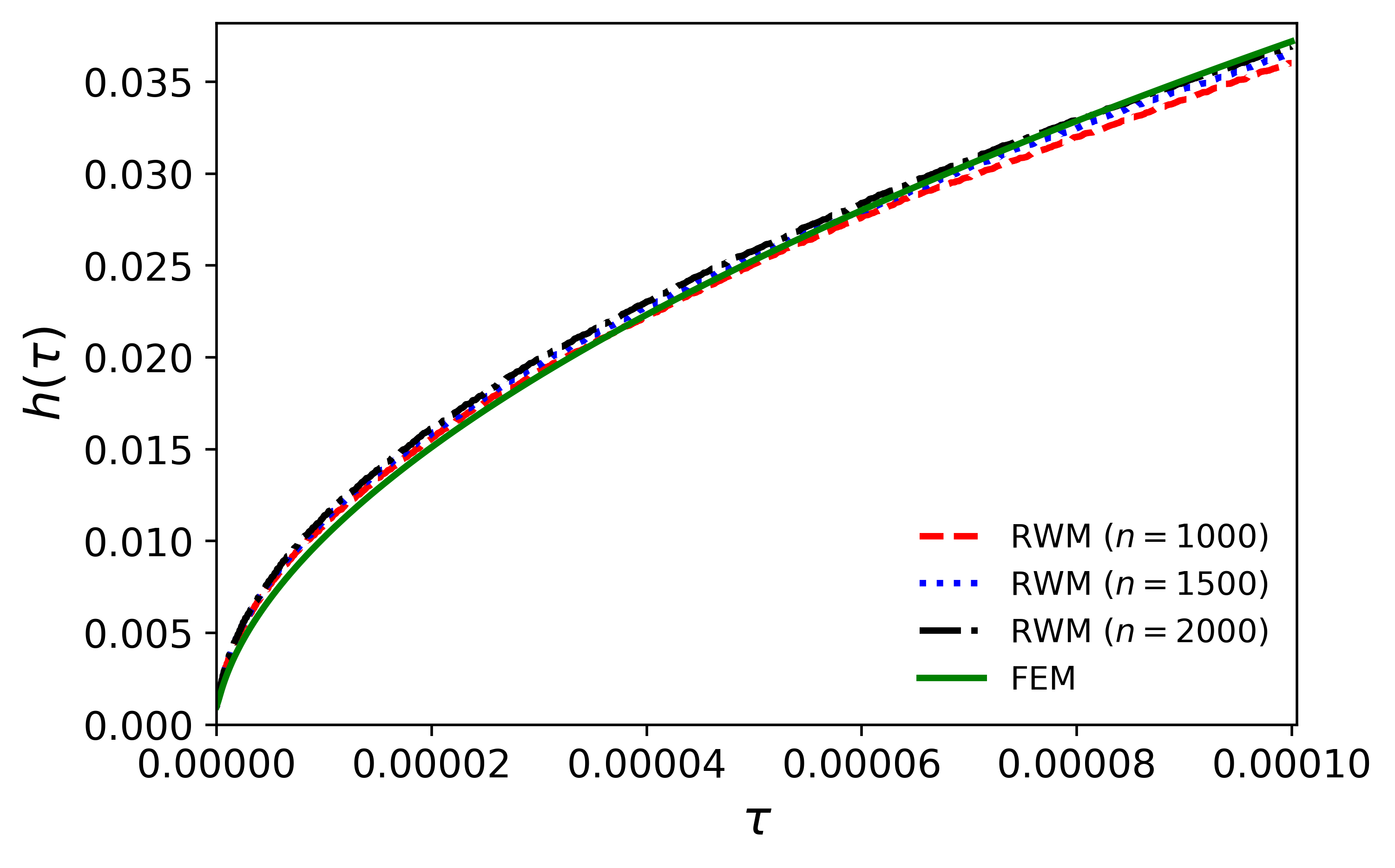}
	\caption{Comparison of the moving front  by RWM for different values of $n$ and  FEM  with   $\Delta \tau = 5 \times 10^{-8}$ (left) and $\Delta \tau = 2.5 \times 10^{-8}$ (right).}
	\label{compare23}
\end{figure}

In Figure \ref{compare22}, we present the concentration profile at time $\tau= 0.00005$  for different values of $n$ and compare them to the finite element solution. The plots in Figure \ref{compare22} show that the random walk solution approaches the finite element solution as $n$ increases. Additionally, increasing the value of $n$ leads to a better approximation for the position of the right (moving) boundary, see Figure~\ref{compare23}.
In Figure~\ref{compare23}, we present the corresponding simulation results for the  position of the moving boundary for different values for $n$. As a clear trend, we observe that the quality of the approximation improves with increasing $n$ and decreasing $\Delta \tau$.

\begin{figure}[ht!] 
\centering 
\includegraphics[width=0.45\textwidth]{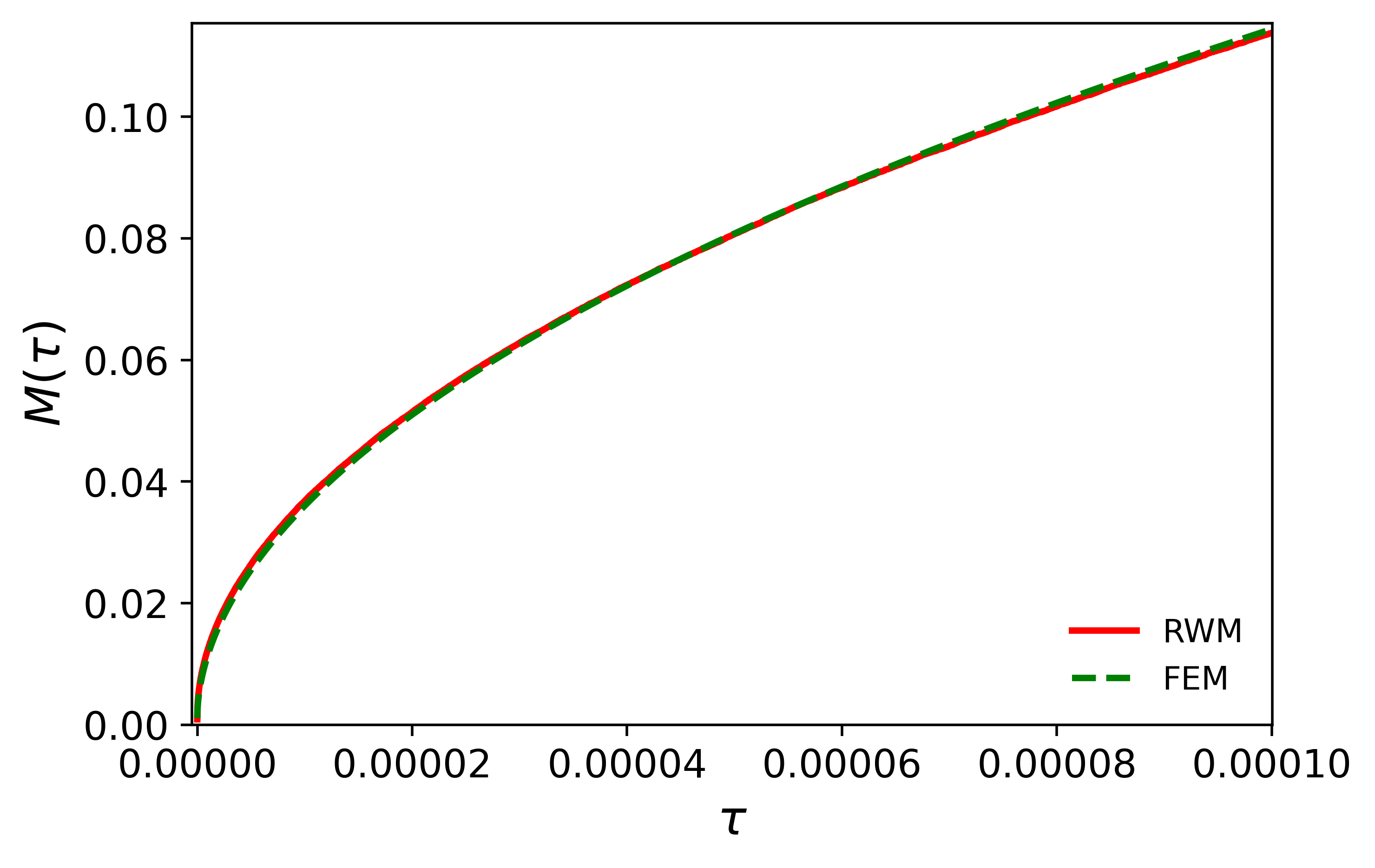}
	\caption{ Comparison of the total mass over time by RWM and FEM, with  $n=2000$, $\Delta \tau = 5 \times 10^{-8}$.}
	\label{compare24}
\end{figure}

We calculate the total mass of the concentration $M(\tau)$ in the domain $(0, h(\tau))$ by 
 \begin{align}
M(\tau) = \int_0^{h(\tau)} u(\tau, z) dz.
\end{align}
 Figure \ref{compare24} shows  the evolution in time of the total mass
for $n = 2000$ and $\Delta \tau = 5 \times 10^{-8}$. Comparing the random walk solution to the finite element solution,  it can be seen that the total mass of concentration obtained by  the finite element method is reproduced by the random walk solution.

 \subsection{Simulation results for Robin boundary condition}
\begin{figure}[ht] 
	\centering
 \hspace{0.01cm}
\includegraphics[width=0.45\textwidth]{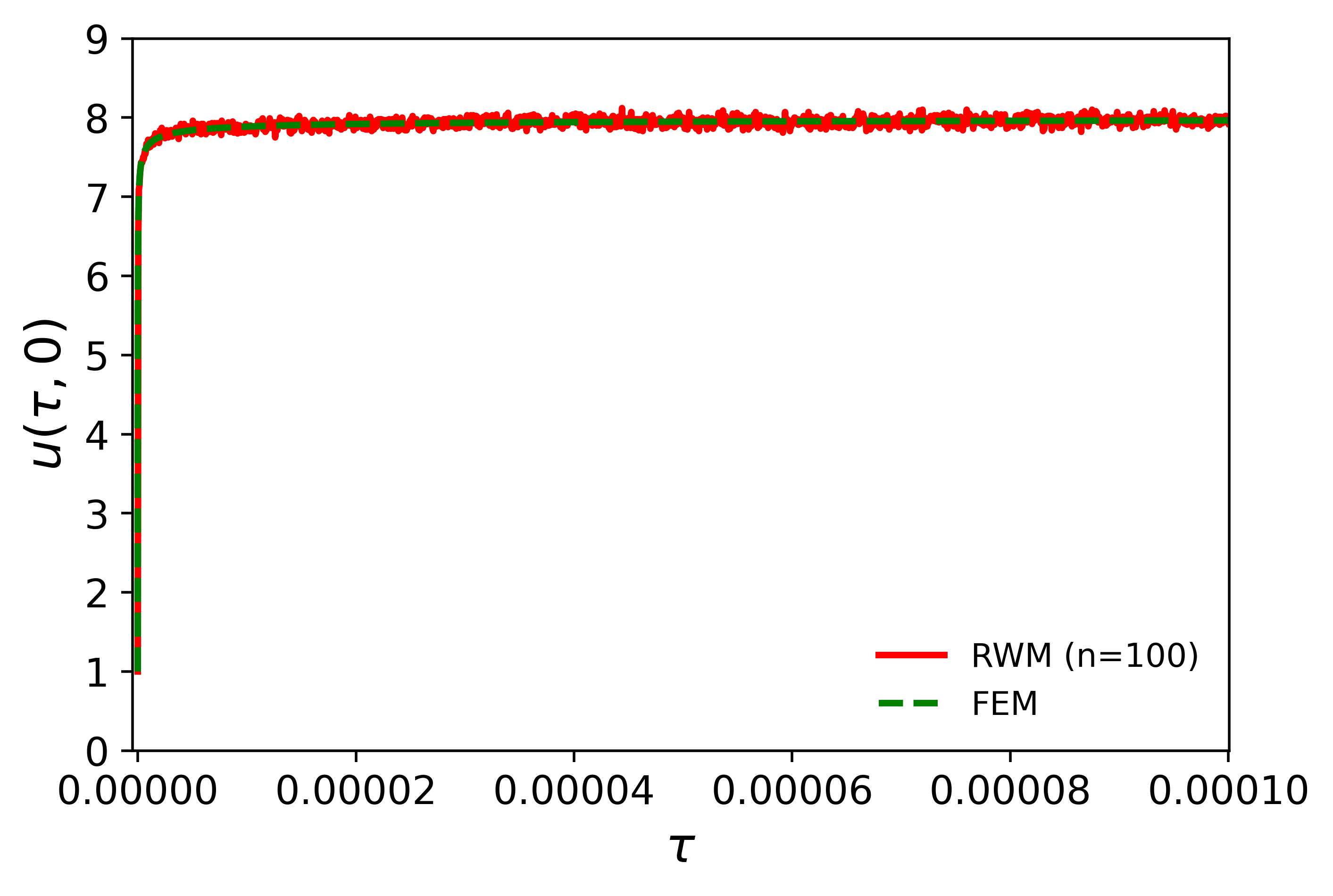}
	\hspace{0.01cm}
\includegraphics[width=0.45\textwidth]
{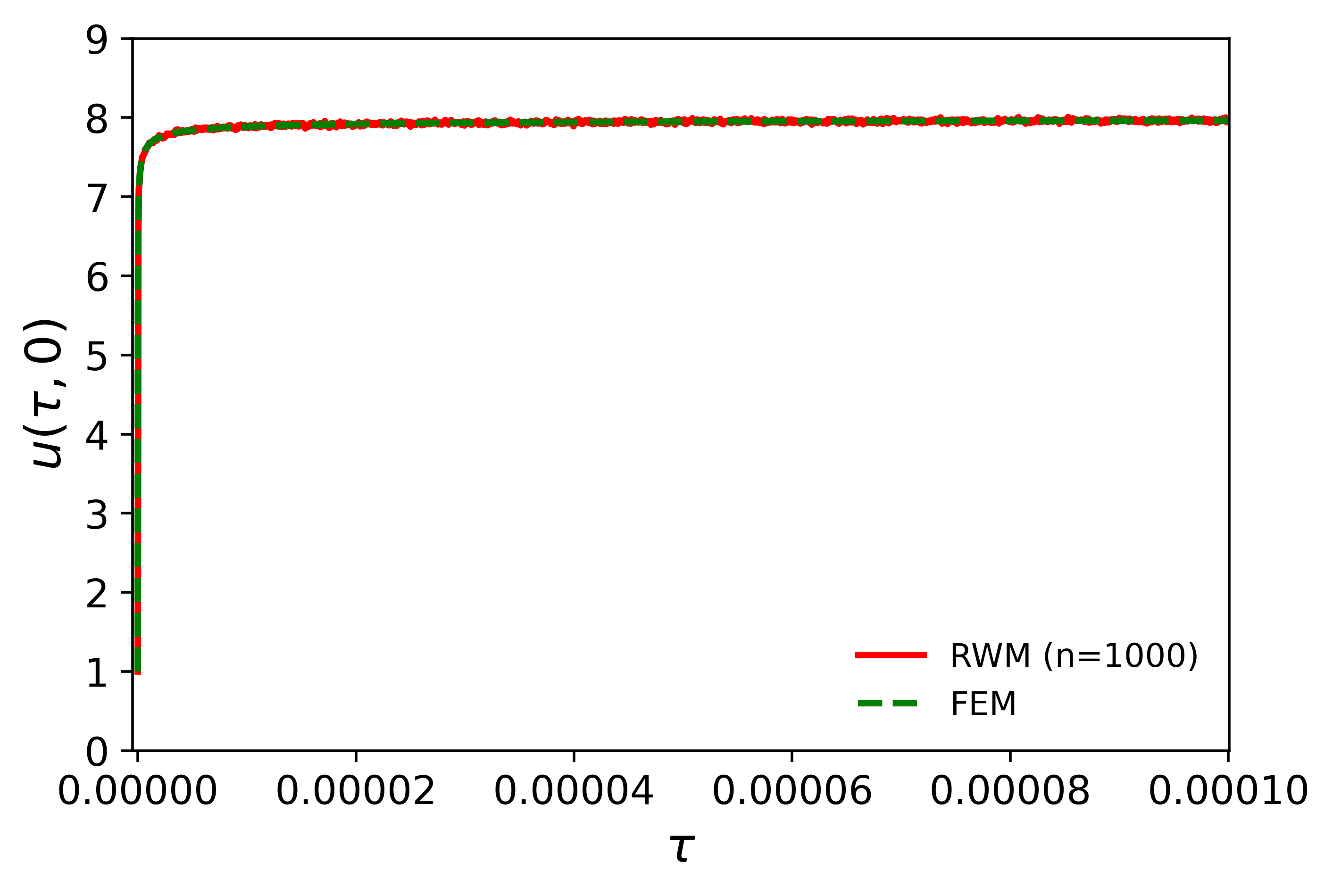}
 \caption{Numerical approximations of the concentration profile at the left boundary $z=0$, for $n=100$ (left), and $n=1000$ (right), with $\Delta \tau = 5 \times 10^{-8}$. }
  \label{randomfemmvleftbc}
\end{figure}

In this section, we present the random walk solution to  \eqref{4a17}--\eqref{4a22} and compared it with the finite element solution.  
We again take the observation time $T = 0.0001$.
All  parameters are taken to be the same as in Section \ref{nonhomdirichlet} except for the parameters arising at the left boundary. 
We now take $5000$ for  Bi, $2.5$ for H and $10$  for $b(\tau)$.
We choose  $u_0 = 1$ so  that  initially, we place $n$ walkers at each discrete node at $z = z_i$ for $ i \in \{0, 1, \cdots, k_0\}$. For the next  time step, we introduce the walkers at the left boundary $z= 0$,  given by \eqref{4a28}. It is worth noting that in Section~\ref{nonhomdirichlet} we prescribed an equal number of walkers at $z=0$ for each time step.  However, the  number of walkers  at $z=0$ now depends on  time. Therefore, to get a better approximation, it is crucial to  introduce correctly the number of the walkers at $z= 0$.  
In Figure \ref{randomfemmvleftbc}, we compare the concentration profile at the left boundary $z=0$ obtained by RWM and FEM. 
The plots in Figure \ref{randomfemmvleftbc} illustrate that with increasing  $n$ the approximations by the two methods are in good agreement.

\begin{figure}[ht] 
	\centering
 \hspace{0.01cm}
\includegraphics[width=0.45\textwidth]{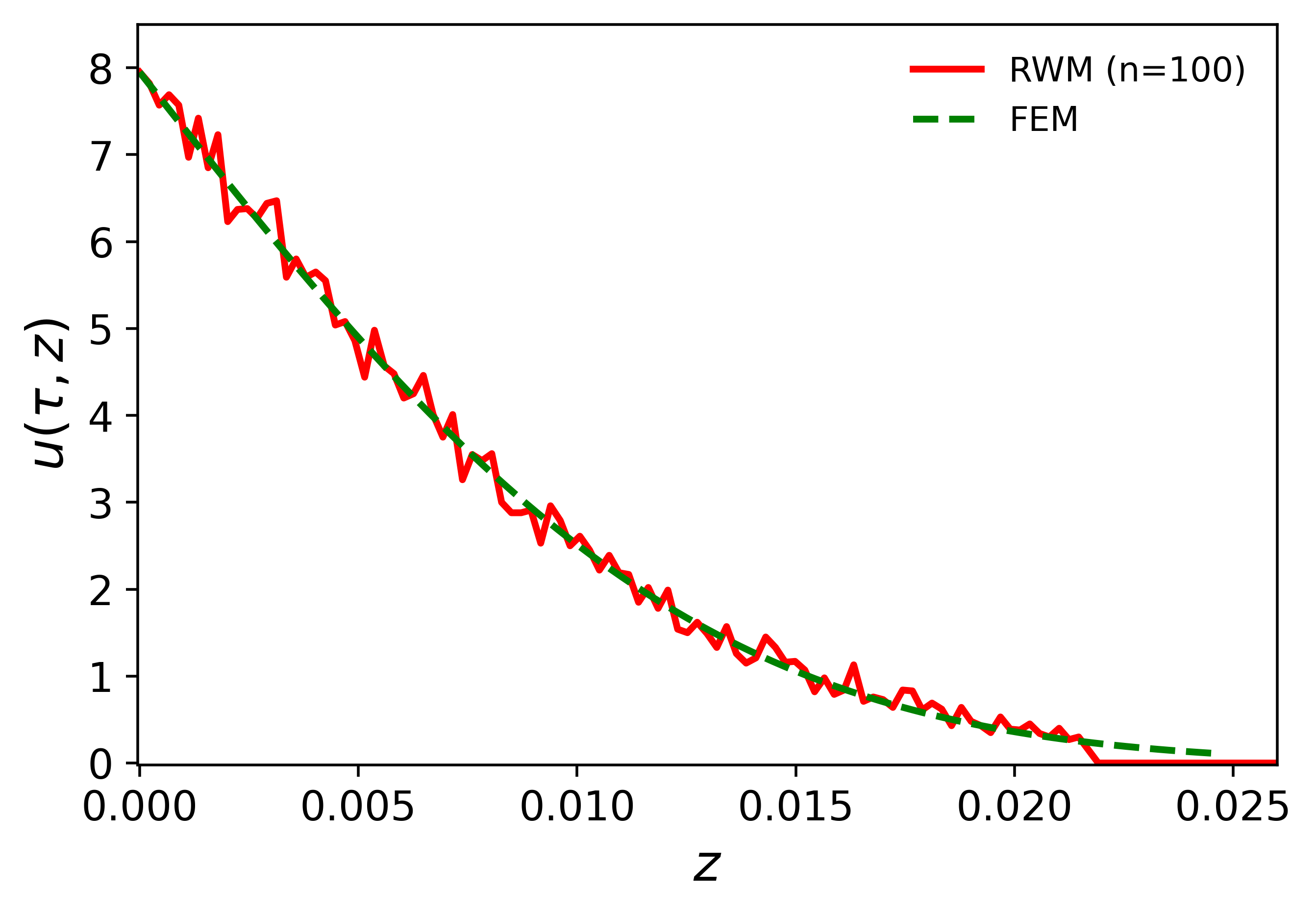}
	%\hspace{0.01cm}
\includegraphics[width=0.45\textwidth]
{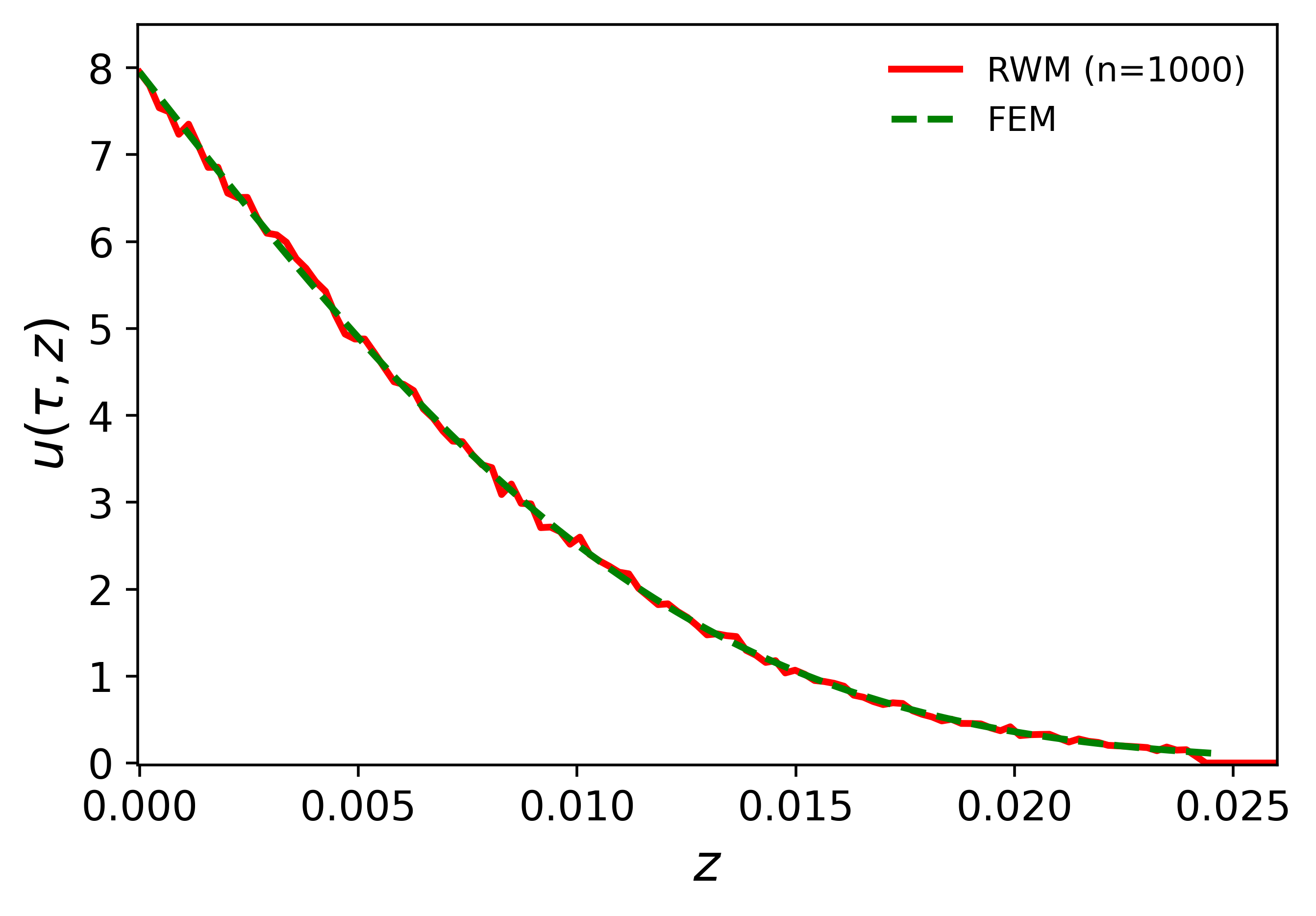}
 \caption{Numerical approximations of the concentration profile at $\tau=0.00005$, for $n=100$ (left), and $n=1000$ (right), with  $\Delta \tau = 2.5 \times 10^{-8}$.} 
 \label{randomfemmvcon}
\end{figure}

In Figure \ref{randomfemmvcon}, we show the concentration profile at $\tau = 0.00005$ obtained by RWM for different values of $n$ and compare them to the solution obtained by FEM. Comparing the plots in Figure \ref{randomfemmvcon} and Figure \ref{compare22}, we observe that the concentration profiles, for the non-homogeneous Dirichlet and Robin-boundary conditions,  have a similar shape.

\begin{figure}[ht!] 
\centering 
\includegraphics[width=0.45\textwidth]{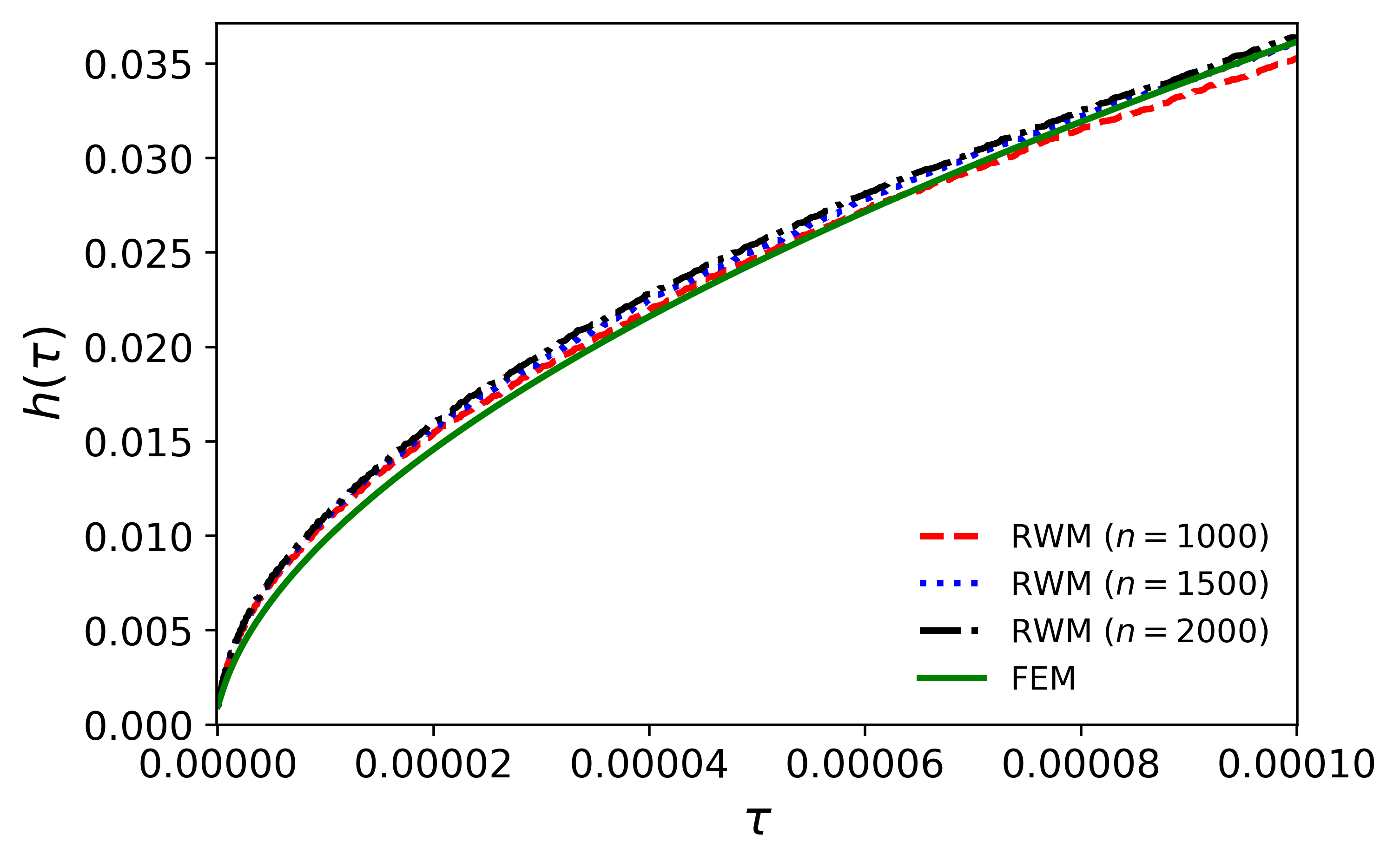}
\hspace{0.001cm}
\includegraphics[width = 0.45\textwidth]{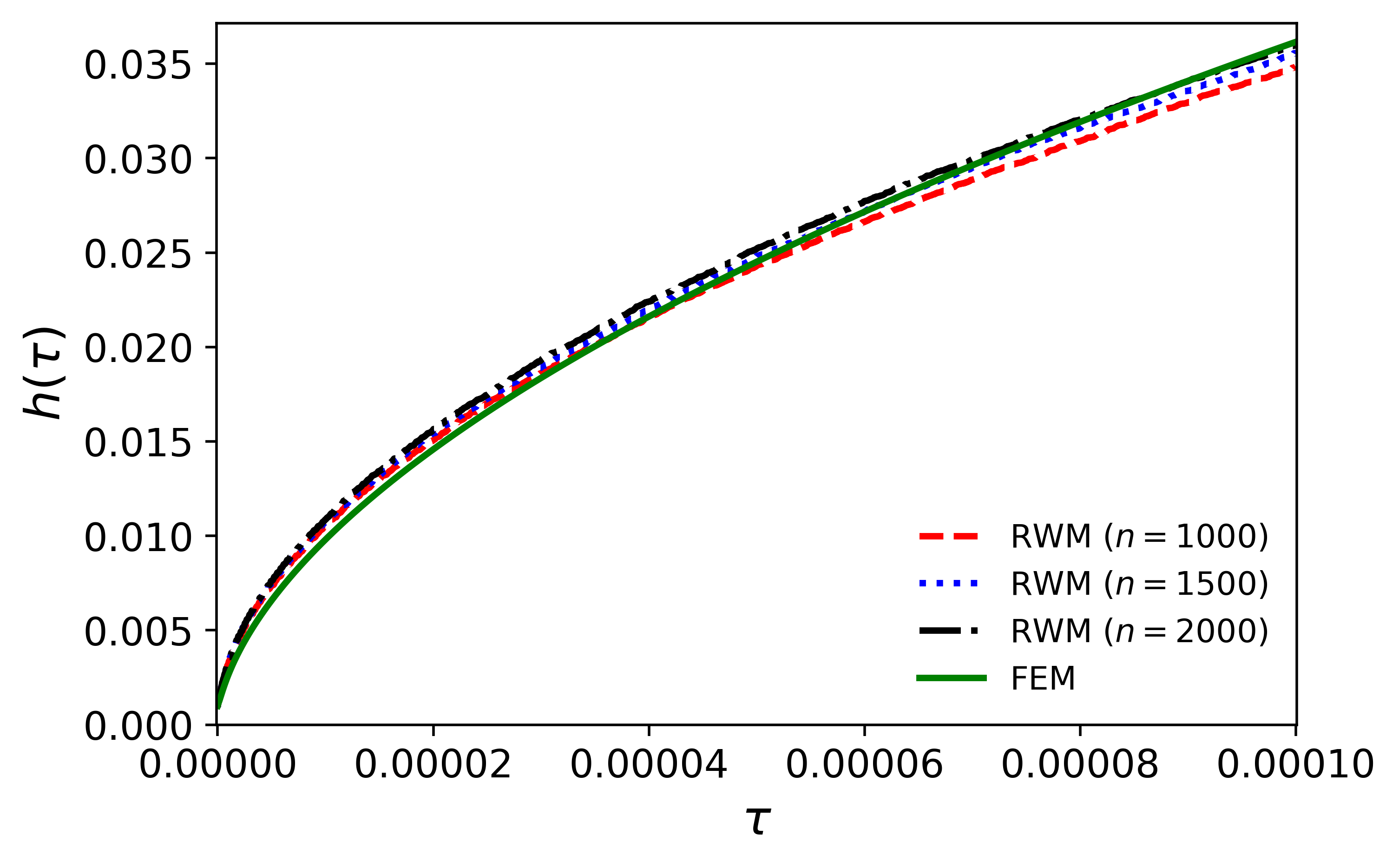}
	\caption{Comparison of the moving front  by RWM for different values of $n$ and  FEM, with $\Delta \tau = 5 \times 10^{-8}$ (left), and $\Delta \tau = 2.5 \times 10^{-8}$ (right).}
	\label{robinbc_diff_n}
\end{figure}

In Figure \ref{robinbc_diff_n}, we present the moving front obtained by RWM for different values of $n$ and compare them with the FEM solution.  
We compute the finite element solution on a time mesh size $\Delta \tau = 5 \times 10^{-9}$ by taking $100$ finite elements.  During the random walk simulation, we observed that as long as the time step $\Delta \tau$ satisfies the requirement \eqref{4a31}, we can ensure that $P_b(\tau_j) < 1$. However, if the walkers at the moving boundary do not satisfy $0<P_b(\tau_j) < 1$, the walkers  contribute neither to the increment of  the moving boundary nor to the increment of the concentration profile at the neighborhood discrete node. If the total number of walkers who do not satisfy $0<P_b(\tau_j) < 1$ is negligibly small, then the contribution of those walkers  does not affect the results of the solution. 
The simulations were done on
%The simulation runs were performed  on 
a Mac with  8 cores 
and 16 GB RAM, using a code written in  \texttt{Python}.

	\begin {table}[ht]
	\begin{center}
		\begin{tabular}{ |p{1.5cm}|p{1.8cm}|p{1.8cm}|p{1.8cm}| p{1.8cm}|}
			\hline
			\diagbox{$\Delta \tau$}{$n$} & $500$ &$1000$&$1500$& $2000$\\
			\hline
			$1.0 \times 10^{-7}$ &1.1139& 2.2317& 3.2732 & 4.4231 \\
   $5.0 \times 10^{-8}$ &3.1012& 6.1230& 9.4097 & 12.2056 \\
		$2.5\times 10^{-8}$ & 8.7836  &17.3675 & 28.3242 & 37.1905\\
			\hline
		\end{tabular}
		\caption {Computational time in minutes for different values of $n$ and $\Delta \tau$.}
		\label{tab:1} 
	\end{center}
	\end {table}

\begin{figure}[ht] 
	\centering
 \hspace{0.01cm}
\includegraphics[width=0.45\textwidth]{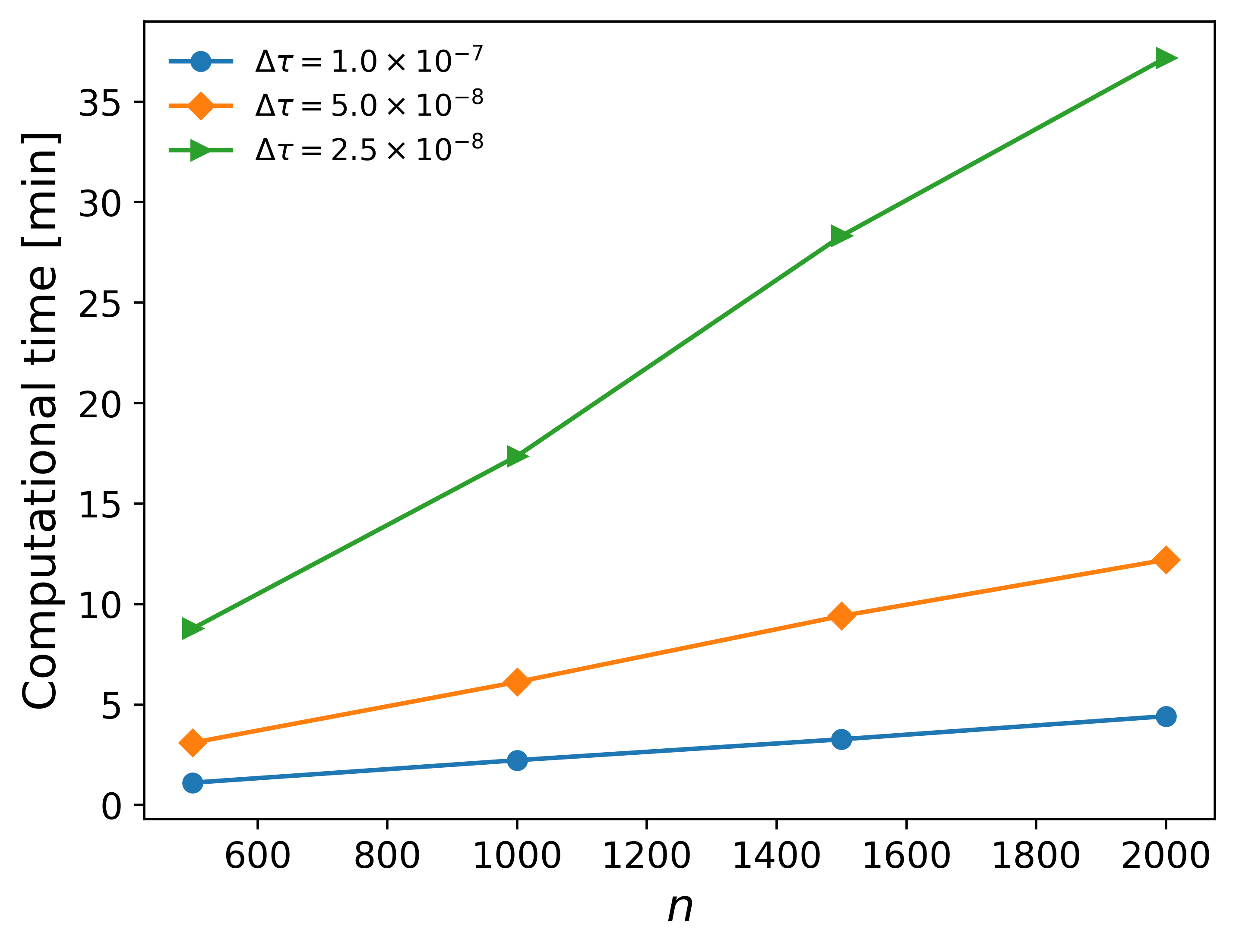}
 \caption{Computational time of RWM  for different values of $n$ and $\Delta \tau$.}
 \label{com_time}
\end{figure}

The computational cost of the simulation mainly depends on  the value of $n$ and on the choice of the time-step size $\Delta \tau$. The effect of time discretization 
 and $n$ is tested by running the simulation for different values of  $\Delta \tau$ and  $n$.  We list the computational time for different values of $n$ and $\Delta \tau$ in Table \ref{tab:1}. To  visualize the dependency of computational time over the chosen $n$,  we plot the corresponding  computational time over $n$  in Figure \ref{com_time}. We observe that  the computational time increases linearly with increasing $n$.

\subsection{Simulation results for capturing  laboratory-controlled penetration depths}
\begin{figure}[ht] 
	\centering
 \hspace{0.01cm}
\includegraphics[width=0.45\textwidth]{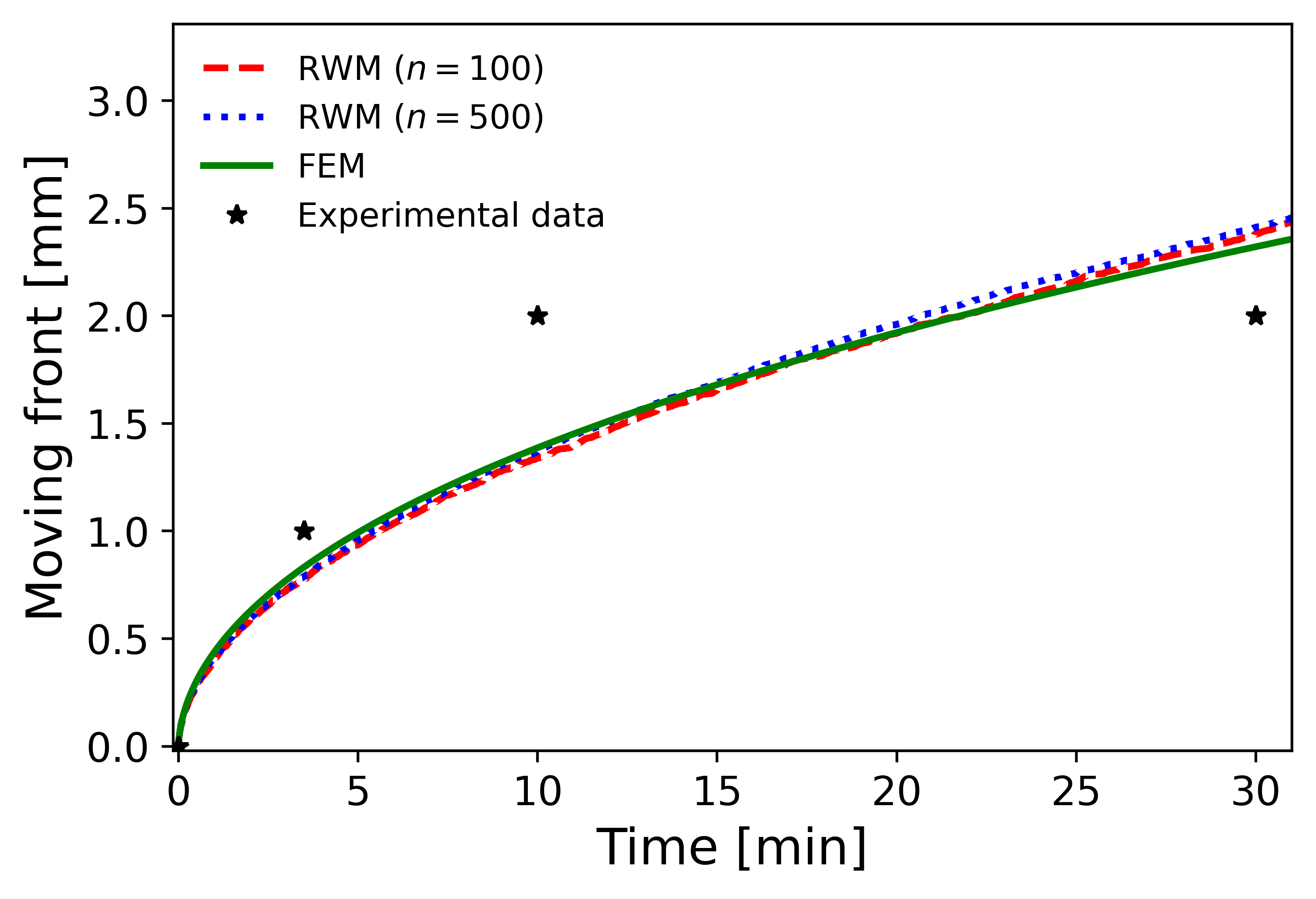}
 \caption{Comparison of the experimental data {\em versus} the computed position of  the moving boundary by both FEM and RWM.}
 \label{moving_lab_fem_rwm}
\end{figure}

%Our next step involves doing  simulations to compare the experimentally measured  positions of the moving boundary against the  position of the moving boundary obtained by FEM and RWM.
In this section, we compare the experimentally measured positions of the moving boundary with the results of our simulations based on the RWM and FEM.
\begin {table}[ht]
\begin{center}
	\begin{tabular}{ |p{7.5cm}|p{1.8cm}|p{3.5cm}| }
		\hline
		Parameters & Dimension&Typical Values\\
		\hline
		Diffusion constant, $D$ &$L^2T^{-1}$& $0.01$ (mm$^2$/min) \\\hline
		Absorption rate,  $\beta$ &$ LT^{-1}$ & $0.564$ (mm/min)\\\hline
		Constant  $a_0$ & $L^4 T^{-1} M^{-1}$& $50$ (mm$^4$/min/gram) \\\hline
		Initial height of diffusants,  $s_0$ & $L$ & $0.01$ (mm)\\\hline
		$\sigma(s(t))$& $ML^{-3}$& $0.5s(t)$ (gram/mm$^3$)\\\hline
		Initial diffusant concentration, $m_0$ & $ML^{-3}$& $0.5$ (gram/mm$^3$)\\\hline
		Concentration in lower surface of the rubber, $b$& $ML^{-3}$&$10$ (gram/mm$^3$) \\\hline
		Henry's constant, 	$\rm{H}$& --& 2.50 (dimensionless)\\
		\hline
	\end{tabular}
	\caption {Name,  dimension and typical values for the reference model  parameters.}
	\label{parameter} 
\end{center}
\end {table}
The typical values of our reference parameters and the dimension of the parameters are listed in Table \ref{parameter}. If not stated otherwise, the numerical results included in this section are produced using the values of the reference parameters. The reference set of parameters for the current simulations is not precisely the same as in our previous work \cite{nepal2021moving}.  Some of them are chosen based on the experimental setup,  while some of them are chosen to capture the laboratory penetration front. We perform our simulations  to recover experimental findings related to the diffusion of cyclohexane and the resulting swelling in a piece of material made of ethylene propylene diene monomer rubber (EPDM). 
We take $T_f = 31$ minutes for the final time. We take the value $10$ mm and $0.5$ gram/mm$^3$ for the characteristic length scale $x_{ref}$ and the reference value for the concentration profile $m_{ref}$, respectively. 
With this choice of parameters, the dimensional numbers Bi and $A_0$ defined in Section \ref{modelequation} are of order $10^4$ and $10^2$, respectively. 
%$O(10^4)$ and $O(10^2)$, respectively.\\
Initially, the diffusant is uniformly distributed within the rubber up to $0.01$ mm. Within a short time, the diffusant 
%quickly 
enters the rubber from the left boundary $x = 0$ and  diffuses further into the material. 
In Figure \ref{moving_lab_fem_rwm}, we present the approximated positions of  the moving boundary computed by both  random walk  and finite element methods.  We compare these computed profiles of penetration depths {\em versus} time against  the experimental data reported in \cite{nepal2021moving}.

\begin{figure}[ht] 
\centering 
\includegraphics[width=0.45\textwidth]{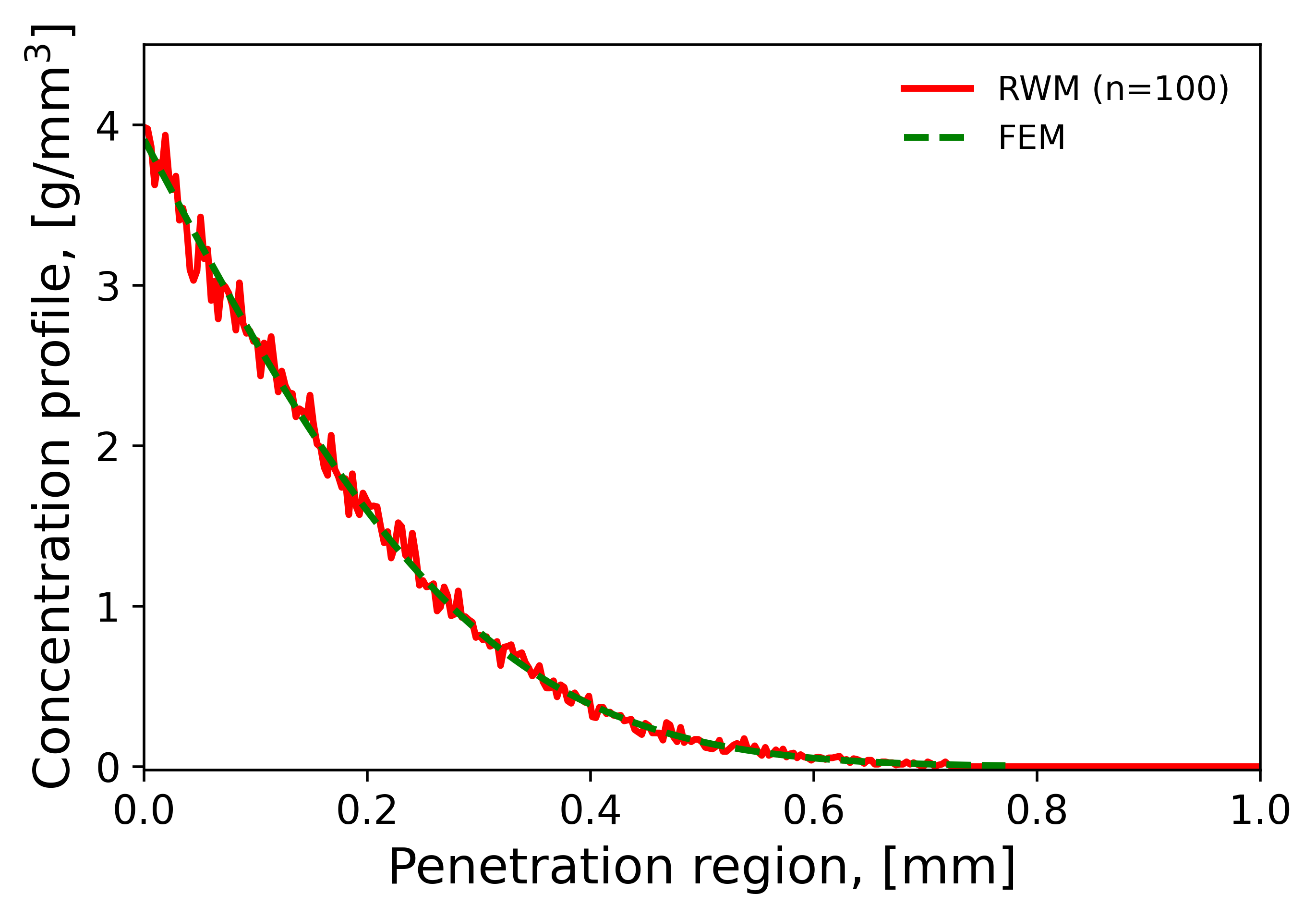}
\hspace{0.01cm}
\includegraphics[width = 0.45\textwidth]{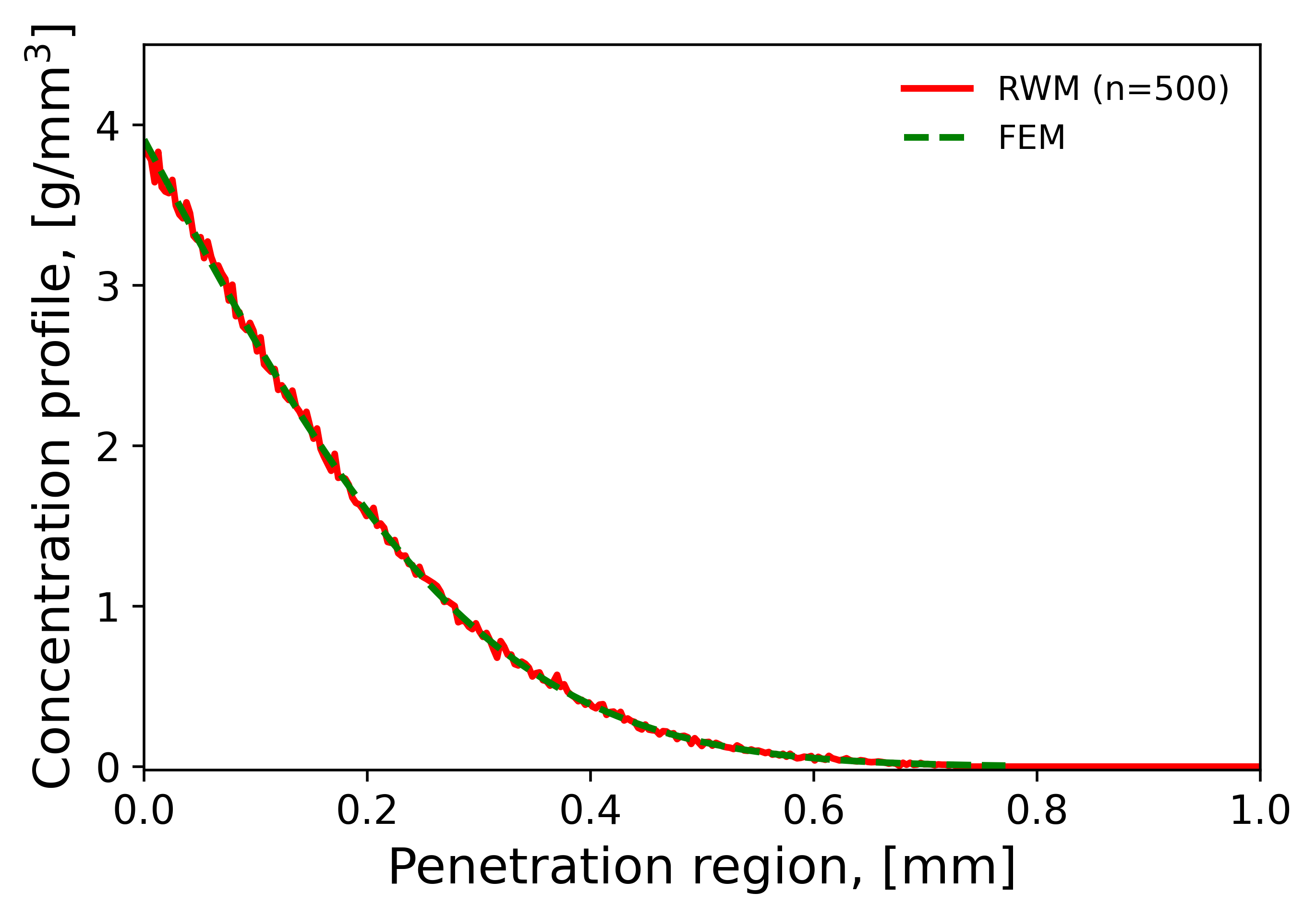}
	\caption{Numerical approximations of the concentration profile at $t = 3$ min, for $n=100$ (left), and $n=500$ (right), with $\Delta t = 0.0005$ min.}
\label{compare_con_fem_rwm}
\end{figure}

In Figure \ref{compare_con_fem_rwm}, we show the concentration profile at $t = 3$ min obtained by RWM for different values of $n$ and compare them to the solution obtained by FEM.

\section{Conclusion and outlook}\label{conclusion}

We reported a random walk method (RWM) capable to approximate, with controllable accuracy, the FEM approximation of the weak solution to our moving boundary problem with a kinetic condition describing the penetration of diffusants into rubber. Since in our previous work \cite{nepal2023analysis}, we did prove convergence rates (in suitable  norms) of the FEM approximation to the original solution, the quality of the FEM approximation transfers to the approximation by RWM.

Trusting \cite{lu1998convergence, hald1981convergence,kaushansky2023convergence}, it seems to be possible to analyse rigorously the convergence property of the proposed RWM to approximate the targeted moving boundary problem. In \cite{lu1998convergence, hald1981convergence}, the  authors solve a reaction-diffusion equation  by a fractional step method, which acts like an operator splitting between computing 
a deterministic ordinary differential equation and running a random walk method for the diffusion equation. 
They prove that the expected value of the computed solution tends to the finite difference approximation of the  (strong) solution to the reaction-diffusion
equation with a suitable control on the variance. 
Additionally, we refer the reader also to the slightly more sophisticated approach mentioned in \cite{kaushansky2023convergence}, which fits very well to our setting when $\sigma(r)=0$ for all $r\geq 0$. 
The main advantage of \cite{kaushansky2023convergence}, compared to e.g. \cite{lu1998convergence}, is that in the first work one points out that, at least in 1D, moving boundary problems with kinetic conditions may admit probabilistic interpretations, which is in principle not {\em a priori} obvious due to the non-dissipative feature of the kinetic condition. As further work, careful mathematical analysis needs to be performed to shed light on an eventual probabilistic interpretation of our problem as described in Section \ref{modelequation}. 

Last but not least, the proposed RWM is able to compute the large time behavior of the moving boundary so that the experimental results reported in \cite{nepal2021moving} are qualitatively recovered. 
The attention here was focused exclusively on the dense rubber case.

\section*{Acknowledgments}
 The activity of S.N.  and A.M. is financed partially by the Swedish Research Council's project {\em "Homogenization and dimension reduction of thin heterogeneous layers"}, grant nr. VR 2018-03648. A.M. also thanks the Knowledge Foundation for the grant KK 2019-0213, which led to the formulation of this problem setting. The authors thank Dr. N. Suciu (Cluj, Romania) for inspiring discussions on closely related matters.
\appendix

 %\newpage
\begin{center}
	\bibliographystyle{plain}
	\bibliography{Random_walk_method}

\begin{thebibliography}{10}

\bibitem{ZAMM}
T.~Aiki, K.~Kumazaki, and A.~Muntean.
\newblock A free boundary problem describing migration into rubbers – {Q}uest
  for the large time behavior.
\newblock {\em Z. Angew. Math. Mech. (ZAMM)}, 102:e202100134, 2022.

\bibitem{boccardo2018improved}
G.~Boccardo, I.M. Sokolov, and A.~Paster.
\newblock An improved scheme for a {R}obin boundary condition in discrete-time
  random walk algorithms.
\newblock {\em Journal of Computational Physics}, 374:1152--1165, 2018.

\bibitem{borges2021effect}
C.S.P. Borges, A.~Akhavan-Safar, E.A.S. Marques, R.J.C. Carbas, C.~Ueffing,
  P.~Wei{\ss}graeber, and L.F.M. da~Silva.
\newblock Effect of water ingress on the mechanical and chemical properties of
  polybutylene terephthalate reinforced with glass fibers.
\newblock {\em Materials}, 14(5):1261, 2021.

\bibitem{CASABAN2022}
M.-C. Casabán, R.~Company, and L.~Jódar.
\newblock Numerical difference solution of moving boundary random {S}tefan
  problems.
\newblock {\em Mathematics and Computers in Simulation}, 2022.

\bibitem{chester2015finite}
S.~A. Chester, C.~V. {Di Leo}, and L.~Anand.
\newblock A finite element implementation of a coupled diffusion-deformation
  theory for elastomeric gels.
\newblock {\em International Journal of Solids and Structures}, 52:1--18, 2015.

\bibitem{Erban_2007}
R.~Erban and S.~J. Chapman.
\newblock Reactive boundary conditions for stochastic simulations of
  reaction–diffusion processes.
\newblock {\em Physical Biology}, 4(1):16 -- 28, 2007.

\bibitem{haji1967solution}
A.~Haji-Sheikh and E.~M. Sparrow.
\newblock The solution of heat conduction problems by probability methods.
\newblock {\em Journal of Heat Transfer}, 89(2):121--130, 1967.

\bibitem{hald1981convergence}
O.~H. Hald.
\newblock Convergence of random methods for a reaction-diffusion equation.
\newblock {\em SIAM Journal on Scientific and Statistical Computing},
  2(1):85--94, 1981.

\bibitem{kaushansky2023convergence}
V.~Kaushansky, C.~Reisinger, M.~Shkolnikov, and Z.~Q. Song.
\newblock Convergence of a time-stepping scheme to the free boundary in the
  supercooled {S}tefan problem.
\newblock {\em The Annals of Applied Probability}, 33(1):274--298, 2023.

\bibitem{NHM}
K.~Kumazaki and A.~Muntean.
\newblock Local weak solvability of a moving boundary problem describing
  swelling along a halfline.
\newblock {\em Networks \& Heterogeneous Media}, 14:445--469, 2019.

\bibitem{kumazaki2020global}
K.~Kumazaki and A.~Muntean.
\newblock Global weak solvability, continuous dependence on data, and large
  time growth of swelling moving interfaces.
\newblock {\em Interfaces and Free Boundaries}, 22(1):27--50, 2020.

\bibitem{lu1998convergence}
W.~Lu.
\newblock Convergence of a random walk method for a partial differential
  equation.
\newblock {\em Mathematics of Computation}, 67(222):593--602, 1998.

\bibitem{nepal2021moving}
S.~Nepal, R.~Meyer, N.~H. Kr{\"o}ger, T.~Aiki, A.~Muntean, Y.~Wondmagegne, and
  U.~Giese.
\newblock A moving boundary approach of capturing diffusants penetration into
  rubber: {FEM} approximation and comparison with laboratory measurements.
\newblock {\em Kautschuk Gummi Kunststoffe}, 5:61--69, 2021.

\bibitem{nepal2021error}
S.~Nepal, Y.~Wondmagegne, and A.~Muntean.
\newblock Error estimates for semi-discrete finite element approximations for a
  moving boundary problem capturing the penetration of diffusants into rubber.
\newblock {\em International Journal of Numerical Analysis \& Modeling},
  19(1):101--125, 2022.

\bibitem{nepal2023analysis}
S.~Nepal, Y.~Wondmagegne, and A.~Muntean.
\newblock Analysis of a fully discrete approximation to a moving-boundary
  problem describing rubber exposed to diffusants.
\newblock {\em Applied Mathematics and Computation}, 442:127733, 2023.
\newblock \url{https://doi.org/10.1016/j.amc.2022.127733}.

\bibitem{Hogren2014local}
M.~{\"O}gren.
\newblock Local boundary conditions for {NMR}-relaxation in digitized porous
  media.
\newblock {\em The European Physical Journal B}, 87(11):1--6, 2014.

\bibitem{ogren2020stochastic}
M.~{\"O}gren.
\newblock Stochastic solutions of {S}tefan problems with general time-dependent
  boundary conditions.
\newblock \url{https://doi.org/10.1007/978-3-031-17820-7_29} In: Malyarenko,
  A., Ni, Y., Rančić, M., Silvestrov, S. (eds) Stochastic Processes,
  Statistical Methods, and Engineering Mathematics. SPAS 2019. Springer
  Proceedings in Mathematics \& Statistics, vol 408. Springer, Cham.
  (arXiv:2006.04939), 2022.

\bibitem{ogren2019numerical}
M.~{\"O}gren, D.~Jha, S.~Dobbersch{\"u}tz, D.~M{\"u}ter, M.~Carlsson,
  M.~Gulliksson, S.L.S. Stipp, and H.O. S{\o}rensen.
\newblock Numerical simulations of {NMR} relaxation in chalk using local
  {R}obin boundary conditions.
\newblock {\em Journal of Magnetic Resonance}, 308:106597, 2019.

\bibitem{rostami2021chemistry}
E.~Rostami-Tapeh-Esmaeil, A.~Vahidifar, E.~Esmizadeh, and D.~Rodrigue.
\newblock Chemistry, processing, properties, and applications of rubber foams.
\newblock {\em Polymers}, 13(10):1565, 2021.

\bibitem{salsa2016partial}
S.~Salsa.
\newblock {\em Partial {D}ifferential {E}quations in {A}ction: {F}rom
  {M}odelling to {T}heory}, volume~99.
\newblock Springer, 2016.

\bibitem{schwind2003some}
M.~Schwind.
\newblock Some remarks on boundary conditions for random walk -- the
  {S}{\"o}derholm condition.
\newblock {\em Scripta Materialia}, 48(4):461--465, 2003.

\bibitem{schwind2001random}
M.~Schwind and J.~{\AA}gren.
\newblock A random walk approach to {O}stwald ripening.
\newblock {\em Acta Materialia}, 49(18):3821--3828, 2001.

\bibitem{suciu2021global}
N.~Suciu, D.~Illiano, A.~Prechtel, and F.~A. Radu.
\newblock Global random walk solvers for fully coupled flow and transport in
  saturated/unsaturated porous media.
\newblock {\em Advances in Water Resources}, 152:103935, 2021.

\bibitem{talebi2017study}
S.~Talebi, K.~Gharehbash, and H.~R. Jalali.
\newblock Study on random walk and its application to solution of heat
  conduction equation by {M}onte {C}arlo method.
\newblock {\em Progress in Nuclear Energy}, 96:18--35, 2017.

\bibitem{wilmers2015continuum}
J.~Wilmers and S.~Bargmann.
\newblock A continuum mechanical model for the description of solvent induced
  swelling in polymeric glasses: Thermomechanics coupled with diffusion.
\newblock {\em European Journal of Mechanics-A/Solids}, 53:10--18, 2015.

\bibitem{yasser2023experimental}
N.~Yasser, A.~Abdelrahman, M.~Kohail, and A.~Moustafa.
\newblock Experimental investigation of durability properties of rubberized
  concrete.
\newblock {\em Ain Shams Engineering Journal}, page 102111, 2023.
\newblock \url{https://doi.org/10.1016/j.asej.2022.102111}.

\end{thebibliography}
\end{center} 

\end{document}